\documentclass[11pt]{article}%
\usepackage{amssymb}
\usepackage{amsfonts}
\usepackage{amsmath}
\usepackage{graphicx}%
\providecommand{\U}[1]{\protect\rule{.1in}{.1in}}
\newtheorem{theorem}{Theorem}

\begin{document}

\title{Decay towards the overall-healthy state in SIS epidemics on networks}
\author{Piet Van Mieghem\thanks{ Faculty of Electrical Engineering, Mathematics and
Computer Science, P.O Box 5031, 2600 GA Delft, The Netherlands; \emph{email}:
P.F.A.VanMieghem@tudelft.nl }}
\date{Delft University of Technology\\
8 July 2014}
\maketitle

\begin{abstract}
The decay rate of SIS epidemics on the complete graph $K_{N}$ is computed
analytically, based on a new, algebraic method to compute the second largest
eigenvalue of a stochastic three-diagonal matrix up to arbitrary precision.
The latter problem has been addressed around 1950, mainly via the theory of
orthogonal polynomials and probability theory. The accurate determination of
the second largest eigenvalue, also called the \emph{decay parameter}, has
been an outstanding problem appearing in general birth-death processes and
random walks. Application of our general framework to SIS epidemics shows that
the maximum average lifetime of an SIS epidemics in any network with $N$ nodes
is not larger (but tight for $K_{N}$) than
\[
E\left[  T\right]  \sim\frac{1}{\delta}\frac{\frac{\tau}{\tau_{c}}\sqrt{2\pi}%
}{\left(  \frac{\tau}{\tau_{c}}-1\right)  ^{2}}\frac{\exp\left(  N\left\{
\log\frac{\tau}{\tau_{c}}+\frac{\tau_{c}}{\tau}-1\right\}  \right)  }{\sqrt
{N}}=O\left(  e^{N\ln\frac{\tau}{\tau_{c}}}\right)
\]
for large $N$ and for an effective infection rate $\tau=\frac{\beta}{\delta}$
above the epidemic threshold $\tau_{c}$. Our order estimate of $E\left[
T\right]  $ sharpens the order estimate $E\left[  T\right]  =O\left(
e^{bN^{a}}\right)  $ of Draief and Massouli\'{e} \cite{Draief_Massoulie}.
Combining the lower bound results of Mountford \emph{et al.}
\cite{Mountford2013} and our upper bound, we conclude that for almost all
graphs, the average time to absorption for $\tau>\tau_{c}$ is $E\left[
T\right]  =O\left(  e^{c_{G}N}\right)  $, where $c_{G}>0$ depends on the
topological structure of the graph $G$ and $\tau$.

\end{abstract}

\section{Introduction}

We consider a simple dynamic process, a Susceptible-Infected-Susceptible (SIS)
epidemic, on an undirected and unweighted graph $G$ with $N$ nodes and $L$
links, that can be represented by a $N\times N$ symmetric adjacency matrix
$A$. In a SIS epidemic process, the viral state of a node $i$ at time $t$ is
specified by a Bernoulli random variable $X_{i}\left(  t\right)  \in\{0,1\}$:
$X_{i}\left(  t\right)  =0$ for a healthy, but susceptible node and
$X_{i}\left(  t\right)  =1$ for an infected node. A node $i$ at time $t$ can
be in one of the two states: \emph{infected}, with probability $v_{i}%
(t)=\Pr[X_{i}(t)=1]$ or \emph{healthy}, with probability $1-v_{i}(t)$, but
susceptible to the infection. We assume that the curing process per node $i$
is a Poisson process with rate $\delta$ and that the infection rate per link
is a Poisson process with rate $\beta$. Obviously, only when a node is
infected, it can infect its direct neighbors, that are still healthy. Both the
curing and infection Poisson process are independent. The effective infection
rate is defined by $\tau=\frac{\beta}{\delta}$. This is the general
continuous-time description of the simplest type of a SIS epidemic process on
a network. This SIS process with curing rate $\delta=1$ is sometimes also
called the contact process.

Kermack and McKendrick \cite{Kermack_McKendrick1927}, whose work is nicely
reviewed in \cite{Breda_DiekmannJBD2012}, have already demonstrated in 1927
that epidemics generally, thus also the SIS process in particular, possess
\textquotedblleft threshold behavior\textquotedblright. For effective
infection rates below the epidemic threshold, $\tau<\tau_{c}$, the
SIS-infection on networks dies out exponentially fast
\cite{PVM_nonMarkovianSIS_2013}, while for $\tau>\tau_{c}$, the infection
becomes endemic, which means that a non-zero fraction of the nodes remains
infected for a very long time. The precise definition (for finite $N$) and the
computation of the SIS epidemic threshold is still an active field of research
\cite{PVM_RMP_epidemics2014}, though a sharp lower bound exists for any graph,
$\tau_{c}\geq\frac{1}{\lambda_{1}}$, where $\lambda_{1}$ is the largest
eigenvalue of the adjacency matrix of the network
\cite{PVM_nonMarkovianSIS_2013}.

Besides the epidemic threshold, the Markovian SIS process also possesses an
important second property: an absorbing state equal to the overall-healthy
state in which the virus has been eradicated from the network. Draief and
Massoulli\'{e} \cite{Draief_Massoulie} prove that the time $T$ for the SIS
Markov process to hit the absorbing state when the effective infection rate
$\tau<\tau_{c}$ is, on average, not larger than $E\left[  T\right]  \leq
\frac{\log N+1}{\delta-\beta\lambda_{1}}$. On the other hand, when $\tau
>\tau_{c}$, they show that the average time to absorption grows for large $N$
as
\begin{equation}
E\left[  T\right]  =O\left(  e^{bN^{a}}\right)
\label{average_absorption_time_Draief_Massoullie}%
\end{equation}
for some constants $a,b>0$. Hence, the average \textquotedblleft
lifetime\textquotedblright\ of the epidemic below and above the epidemic phase
transition are hugely different, which is a general characteristic of a phase
transition. Mountford \emph{et al.} \cite{Mountford2013} proved that, above
the epidemic threshold in trees with bounded degree, i.e. the maximum degree
$d_{\max}<a$, where $a$ is finite, but $d_{\max}\geq2$ (thus excluding e.g.
the star), $E\left[  T\right]  =O(e^{cN})$ for large $N$ and a real number
$c>0$. Moreover, improving a result of Chatterjee and Durrett
\cite{Chatterjee_Durret2009}, they show that for any $\tau>0$ and large $N$,
the time to absorption or extinction on a power law graph grows exponentially
in $N$.

Fill \cite{Fill_JTP2009} gave a nice stochastic interpretation of the time $T$
to absorption in a continuous-time birth and death process with an absorbing
state zero and $N$ other states, described by the infinitestimal generator
$Q$. Given that the process starts in state $N$, then the absorption time $T$
is equal to a sum of independent exponential random variables, whose rates are
the nonzero eigenvalues of $-Q$. Miclo \cite{Miclo_ESAIM2010} has extended
Fill's result to a finite Markov chain, which is irreducible and reversible
outside the absorbing point. Very recently, Economou \emph{et al.}
\cite{Economou_PhysicaA2014} have analysed the SIS model with heterogeneous
infection rates via a block matrix formalism. In their analysis, they gave the
general expression for distribution of the absorption time $T$ as $\Pr\left[
T\leq t\right]  =1-\left(  x_{0}^{T}e^{Q^{\ast}t}\right)  _{2^{N}-1}$, where
$x_{0}$ is the column vector with the initial states and $Q^{\ast}$ is the
infinitesimal generator (see \cite{Economou_PhysicaA2014} for the labelling of
states) in which the row and column corresponding to the absorbing state are
removed. Artalejo \cite{Artalejo_PhysicaA2012} has shown that the time $T_{q}$
to extinction from the quasi-stationary (or metastable) state obeys
$\Pr\left[  T_{q}\leq t\right]  =1-e^{\zeta t}$, where $\zeta\leq0$ is the
second largest eigenvalue\footnote{More precisely, the largest eigenvalue of
the submatrix $Q_{S}$ of $Q$ associated to the transient and finte set $S$ of
states, that is assumed to be irreducible. When the latter condition of
irreducibility is omitted, $\Pr\left[  T_{q}\leq t\right]  $ is still
exponentially distributed \cite[Theorem 1]{Artalejo_PhysicaA2012}, but with a
more complicated mean $E\left[  T_{q}\right]  $ than $\frac{1}{\zeta}$.} of
the infinitesimal generator $Q$.

Here, we derive a sharper estimate than $E\left[  T\right]  =O\left(
e^{bN^{a}}\right)  $ for the longest possible mean absorption time in any
graph, by computing the spectral decomposition of a tri-diagonal, stochastic
matrix $P$ in (\ref{P_general_tridiagonal_band_matrix}), which is presented in
Appendix \ref{sec_general_tri-diagonal_matrices}. Invoking the Lagrange series
on the characteristic polynomial of $P$, the second largest eigenvalue
$1+\zeta$ (with $\zeta\leq0$) of $P$ is deduced in Appendix
\ref{sec_properties_zeros_characteristic_coefficient_P}. Generally, for the
state vector $x\left[  k\right]  $ of a discrete-time Markov process at
discrete-time $k$ with a real second largest eigenvalue, it holds that any
vector norm $\left\vert \left\vert x\left[  k\right]  -\pi\right\vert
\right\vert \sim\left(  1+\zeta\right)  ^{k}+O\left(  \left\vert
1+z_{3}\right\vert ^{k}\right)  $, where $\pi$ is the corresponding
steady-state vector and $z_{3}$ is the third largest (in absolute value)
eigenvalue of $P$. The number of infected nodes in a SIS epidemic process on
the complete graph can be determined \cite{PVM_EpsilonSIS_PRE2012} via a birth
and death process, the continuous-time variant of a general random walk, whose
infinitesimal generator $Q$ is a tri-diagonal matrix. As shown in Section
\ref{sec_epsilon_SIS_epidemics}, the second largest eigenvalue $\zeta$ of the
infinitesimal generator $Q$ can thus be interpreted as the decay rate of the
SIS epidemics on the complete graph towards the overall-healthy state and,
approximately, the average lifetime of an SIS epidemics is about $E\left[
T\right]  \simeq\frac{1}{\left\vert \zeta\right\vert }$. Now, given a fixed
infection rate $\beta$ and curing rate $\delta$, among all networks with $N$
nodes, the SIS infection spreads fastest in the complete graph $K_{N}$ with
$N$ nodes, because each node can be infected by a maximum possible number of
neighbors. Hence, the longest time $T$ to hit the overall-healthy state and,
equivalently, the minimum decay rate $\zeta$ among all graphs are attained in
the complete graph $K_{N}$.

Our main result for SIS epidemics is the accurate expression of the decay rate
$\zeta$ in $K_{N}$ for effective infection rates $\tau>\tau_{c}$ and large
$N$
\begin{equation}
-\zeta=\frac{1}{F\left(  \tau\right)  }+O\left(  \frac{N^{2}\log N}{x^{2N-1}%
}\right)  \label{zeta_SIS_KN_above_threshold_x>1}%
\end{equation}
where $x=\tau N\simeq\frac{\tau}{\tau_{c}}>1$ and where%
\begin{equation}
F\left(  \tau\right)  =\frac{1}{\delta}\sum_{j=1}^{N}\sum_{r=0}^{j-1}%
\frac{\left(  N-j+r\right)  !}{j\left(  N-j\right)  !}\tau^{r}\label{f1_eps=0}%
\end{equation}
The double sum in (\ref{f1_eps=0}) is hard to compute for large $N$ and, after
surprisingly much effort as illustrated in Appendix \ref{sec_F(tau)}, we
established in Theorem \ref{theorem_bound_f1_eps=0} the correct behavior of%
\begin{equation}
F\left(  \frac{x}{N}\right)  \sim\frac{1}{\delta}\frac{x\sqrt{2\pi}}{\left(
x-1\right)  ^{2}}\frac{\exp\left(  N\left\{  \log x+\frac{1}{x}-1\right\}
\right)  }{\sqrt{N}}\label{asymptotic_F(tau)_Nlarge_x_fixed}%
\end{equation}
for large $N$ and fixed $x=\tau N>1$. Roughly, for $x$ slightly above than 1,
we deduce from the asymptotic expression
(\ref{asymptotic_F(tau)_Nlarge_x_fixed}) of $F\left(  \tau\right)  $ that
$E\left[  T\right]  =O\left(  e^{N\ln\frac{\tau}{\tau_{c}}}\right)  $. The
exponentially accurate order estimate (\ref{zeta_SIS_KN_above_threshold_x>1})
thus specifies the parameters $a=1$ and $b=\ln\frac{\tau}{\tau_{c}}$ (or more
correctly $b=\ln\frac{\tau}{\tau_{c}}+\frac{\tau_{c}}{\tau}-1$) in the general
estimate (\ref{average_absorption_time_Draief_Massoullie}). Earlier in
\cite{PVM_ToN_VirusSpread}, we have derived the exact $2^{N}\times2^{N}$
infinitesimal generator $Q$ for an SIS process on any graph and have
numerically computed the second smallest eigenvalue of $Q$ for the complete
graph. For small networks up to $N=13$, fitting results suggested that
$E\left[  T\right]  =O\left(  e^{b\left(  \tau\right)  N^{2}}\right)  $.
Hence, the current analytic result $E\left[  T\right]  =O\left(  e^{N\ln
\frac{\tau}{\tau_{c}}}\right)  $ shows that $a\leq1$, in contrast to our
earlier extrapolated order estimates that hinted at $a\leq2$.

The probabilistic interpretation of the absorption time $T$ by Fill
\cite{Fill_JTP2009} leads us to conclude that $E\left[  T\right]  =F\left(
\tau\right)  $, for any value of the effective infection rate $\tau$ (and not,
as above in (\ref{zeta_SIS_KN_above_threshold_x>1}), only for $\tau>\tau_{c}%
$). Thus, starting from the all-infected state, the \emph{exact}\footnote{The
exact relation $E\left[  T\right]  =F\left(  \tau\right)  $ has been verified
by using a hitting time analysis in \cite{PVM_Survival_time_PRE2014}.} average
absorption time $T$ in the SIS process on the complete graph is given by
$F\left(  \tau\right)  $ in (\ref{f1_eps=0}). Moreover, as shown in Appendix
\ref{sec_properties_zeros_characteristic_coefficient_P}, the first term in the
Lagrange series for the second largest eigenvalue $\zeta$ of an infinitesimal
generator $Q$ equals the inverse of the sum of the inverse (non-zero)
eigenvalues of $Q$, which may suggest that, in general Markov chains with an
absorbing state, $-\zeta=\frac{1}{E\left[  T\right]  }+r$, where $r$ are
higher order terms in the Lagrange series.

Finally, combining the lower bound results in \cite{Mountford2013} and the
upper bound in (\ref{zeta_SIS_KN_above_threshold_x>1}), we conclude that for
almost all graphs, the average time to absorption for $\tau>\tau_{c}$ is
$E\left[  T\right]  =O\left(  e^{c_{G}N}\right)  $, where $c_{G}>0$ depends on
the topological structure of the graph $G$. The interesting open next question
lies in the accurate determination of $c_{G}$ for a given graph $G$, different
from $K_{N}$.

\section{Markovian $\varepsilon-$SIS epidemics}

\label{sec_epsilon_SIS_epidemics}We first define the Markovian $\varepsilon
-$SIS epidemics on networks. Besides an infection process with rate $\beta$
per infected neighbor and a nodal curing process with rate $\delta$ as in the
SIS model, each node contains a Poissonean self-infection process with rate
$\varepsilon$. All three Poisson processes are independent. This $\varepsilon
$-SIS epidemic process on the complete graph $K_{N}$ is a birth and death
process with birth rate $\lambda_{j}=\left(  \beta j+\varepsilon\right)
\left(  N-j\right)  $ and death rate $\mu_{j}=j\delta$, as shown in
\cite{PVM_EpsilonSIS_PRE2012}. When the process $X\left(  t\right)  $ at time
$t$ is at state $j$, precisely $j$ nodes in $K_{N}$ are infected. For
$\varepsilon>0$, all rates are positive and the birth and death process is
irreducible, i.e. without absorbing state. Thus, the theory developed in
Appendix \ref{sec_general_tri-diagonal_matrices} is applicable when we
substitute%
\begin{align*}
p_{j} &  \rightarrow\left(  \beta j+\varepsilon\right)  \left(  N-j\right)  \\
q_{j} &  \rightarrow j\delta
\end{align*}
In an irreducible, $n$ states, continuous-time Markov process, the $1\times n$
state vector $s\left(  t\right)  $, with component $i$ equal to $s_{i}\left(
t\right)  =\Pr\left[  X\left(  t\right)  =i\right]  $, satisfies%
\[
s\left(  t\right)  =s\left(  0\right)  e^{Qt}%
\]
where the spectral decomposition of the $n\times n$ matrix (see e.g.
\cite{PVM_PerformanceAnalysisCUP}) is%
\[
e^{Qt}=u\pi+\sum_{j=2}^{n}e^{\mu_{j}t}x_{j}y_{j}^{T}%
\]
and $x_{j}$ and $y_{j}$ are the $n\times1$ right- and left-eigenvector
belonging to the $i$-th largest eigenvalue $\mu_{j}$ of $Q$. The
right-eigenvector belonging to the largest eigenvalue $\mu_{1}=0$ of $Q$ is
$x_{0}=u$, the all-one vector. We denote $\mu_{2}=\zeta$. For large $t$, the
tendency of $s\left(  t\right)  $ towards the steady-state vector $\pi$ equals%
\[
s\left(  t\right)  -\pi\approx we^{\zeta t}%
\]
where $s\left(  0\right)  u=1$ and $w=s\left(  0\right)  x_{2}y_{2}^{T}$ is
not a function of the time $t$.

In Appendix \ref{sec_properties_zeros_characteristic_coefficient_P}, we
demonstrate the general bound $\zeta<-\frac{f_{0}}{f_{1}}$, where $f_{k}$ are
the coefficients (\ref{def_f_k}) of the characteristic polynomial of a
tri-band matrix. Hence, the continuous-time Markov process specified by a
tri-diagonal infinitesimal generator $Q$ converges always faster to the
steady-state than $O\left(  \exp\left(  -\frac{f_{0}}{f_{1}}t\right)  \right)
$. This means \cite{PVM_EpsilonSIS_PRE2012} for an $\varepsilon-$SIS-epidemic
process on the complete graph that the epidemics tends to the SIS
\emph{metastable} state with a time constant faster than $T\left(
\varepsilon\right)  =\frac{f_{1}}{f_{0}}$ time units. In the limit
$\varepsilon\rightarrow0$, where the $\varepsilon-$SIS-epidemic process
behaves as the classical SIS epidemics in which the steady-state is the
overall healthy state (which is the absorbing state for the SIS Markov
process), the decay rate of the epidemics towards this absorbing state is
never slower than $\frac{1}{T\left(  0\right)  }$.

The remainder of this section consists of (a) the determination of the
coefficients $f_{k}$ for the $\varepsilon$-SIS epidemic process on the
complete graph $K_{N}$ (Section \ref{sec_coeff_eps_SIS}), (b) the limit form
of the these coefficients for large $N$ and (c) the three regimes, depending
on whether $\tau\geq\tau_{c}$,$\tau\simeq\tau_{c}$ and $\tau<\tau_{c}$, of the
resulting decay rate $\zeta$ in SIS epidemics ( $\varepsilon\rightarrow0$) for
large $N$, in which our main result (\ref{zeta_SIS_KN_above_threshold_x>1}) is derived.

\subsection{Coefficients $f_{0}$, $f_{1}$ and $f_{2}$ in $\varepsilon-$SIS
epidemics}

\label{sec_coeff_eps_SIS}The inverse of the probability that no node in
$K_{N}$ is infected is \cite{PVM_EpsilonSIS_PRE2012}
\[
f_{0}=\frac{1}{\pi_{0}}=\sum_{k=0}^{N}\binom{N}{k}\tau^{k}\frac{\Gamma\left(
\frac{\varepsilon^{\ast}}{\tau}+k\right)  }{\Gamma\left(  \frac{\varepsilon
^{\ast}}{\tau}\right)  }=1+\sum_{k=1}^{N}\binom{N}{k}\tau^{k}\frac
{\Gamma\left(  \frac{\varepsilon^{\ast}}{\tau}+k\right)  }{\Gamma\left(
\frac{\varepsilon^{\ast}}{\tau}\right)  }%
\]
where $\varepsilon^{\ast}=\frac{\varepsilon}{\delta}$ and $\tau=\frac{\beta
}{\delta}$. Since $\lim_{\varepsilon\rightarrow0}\frac{1}{\Gamma\left(
\frac{\varepsilon^{\ast}}{\tau}\right)  }=0$,
\[
\lim_{\varepsilon\rightarrow0}f_{0}=1
\]
agreeing with the fact that the steady-state in Markovian SIS epidemics is
equal to the overall-healthy state, which is the absorbing state zero. Using%
\[
\prod_{m=0}^{j-1}q_{m+1}=\prod_{m=0}^{j-1}\left(  m+1\right)  \delta
=j!\delta^{j}%
\]
and
\[
\prod_{m=0}^{j-1}p_{m}=\prod_{m=0}^{j-1}\left(  \beta m+\varepsilon\right)
\left(  N-m\right)  =\frac{N!\beta^{j}}{\left(  N-j\right)  !}\prod
_{m=0}^{j-1}\left(  m+\frac{\varepsilon}{\beta}\right)  =\frac{N!\beta
^{j}\Gamma\left(  \frac{\varepsilon}{\beta}+j\right)  }{\left(  N-j\right)
!\Gamma\left(  \frac{\varepsilon}{\beta}\right)  }%
\]
into the general expression (\ref{f1}) for $f_{1}$ yields%
\[
f_{1}=\frac{1}{\delta}\sum_{j=1}^{N}\frac{\tau^{j-1}}{j}\sum_{r=0}^{j-1}%
\sum_{k=0}^{j-1-r}\frac{\binom{N-j+r}{r}\binom{N}{j-1-r-k}}{\binom{j-1}{r}%
}\frac{\Gamma\left(  \frac{\varepsilon^{\ast}}{\tau}+j\right)  \Gamma\left(
\frac{\varepsilon^{\ast}}{\tau}+j-1-r-k\right)  }{\Gamma\left(  \frac
{\varepsilon^{\ast}}{\tau}+j-r\right)  \Gamma\left(  \frac{\varepsilon^{\ast}%
}{\tau}\right)  }\frac{1}{\tau^{k}}%
\]
and $\lim_{\varepsilon\rightarrow0}f_{1}=F\left(  \tau\right)  $, specified in
(\ref{f1_eps=0}).

From the definition (\ref{def_f_k}),
\[
f_{2}=\frac{1}{q_{1}q_{2}}+\sum_{j=3}^{N}\frac{c_{2}(j)}{\prod_{m=0}%
^{j-1}q_{m+1}}%
\]
where $c_{2}\left(  j\right)  $ follows from (\ref{c2(j)_explicit}), we arrive
at%
\begin{align*}
f_{2}  &  =\frac{1}{2\delta^{2}}+\frac{1}{\delta^{2}}\sum_{j=3}^{N}%
\frac{\left(  N-2\right)  !\tau^{j-2}\Gamma\left(  \frac{\varepsilon^{\ast}%
}{\tau}+j\right)  }{j!\left(  N-j\right)  !\Gamma\left(  \frac{\varepsilon
^{\ast}}{\tau}+2\right)  }\\
&  \hspace{0.5cm}+\left(  \frac{\tau\left(  N-1\right)  +\varepsilon^{\ast
}\left(  2N-1\right)  +3}{2\delta^{2}}\right)  \sum_{j=3}^{N}\sum_{k=3}%
^{j}\frac{\left(  k-1\right)  !\left(  N-k\right)  !\tau^{j-k}\Gamma\left(
\frac{\varepsilon^{\ast}}{\tau}+j\right)  }{j!\left(  N-j\right)
!\Gamma\left(  \frac{\varepsilon^{\ast}}{\tau}+k\right)  }\\
&  \hspace{0.5cm}+\frac{1}{\delta^{2}}\sum_{j=3}^{N}\sum_{k=3}^{j}\sum
_{s=1}^{k-3}\sum_{l_{1}=0}^{k-s-1}\sum_{l_{2}=0}^{k-s-l_{1}-1}\frac
{\Gamma\left(  \frac{\varepsilon^{\ast}}{\tau}+k-s-1-l_{1}-l_{2}\right)
\Gamma\left(  \frac{\varepsilon^{\ast}}{\tau}+k-s\right)  \Gamma\left(
\frac{\varepsilon^{\ast}}{\tau}+j\right)  }{\Gamma\left(  \frac{\varepsilon
^{\ast}}{\tau}+k-s-l_{1}\right)  \Gamma\left(  \frac{\varepsilon^{\ast}}{\tau
}\right)  \Gamma\left(  \frac{\varepsilon^{\ast}}{\tau}+k\right)  }\\
&  \hspace{0.5cm}\times\frac{N!\left(  N-\left(  k-s-l_{1}\right)  \right)
!\left(  N-k\right)  !\left(  k-s-1-l_{1}\right)  !\left(  k-1\right)
!\tau^{j-1-s-l_{2}}}{j!\left(  N-\left(  k-s-1-l_{1}-l_{2}\right)  \right)
!\left(  N-\left(  k-s\right)  \right)  !\left(  N-j\right)  !\left(
k-s-l_{1}-l_{2}-1\right)  !\left(  k-s\right)  !}%
\end{align*}
After some tedious calculations, we find that%
\begin{align}
\lim_{\varepsilon\rightarrow0}f_{2}  &  =\frac{1}{2\delta^{2}}+\frac{1}%
{\delta^{2}}\sum_{j=3}^{N}\frac{\left(  N-2\right)  !\tau^{j-2}}{j\left(
N-j\right)  !}+\left(  \frac{\tau\left(  N-1\right)  +3}{2\delta^{2}}\right)
\sum_{j=3}^{N}\sum_{k=3}^{j}\frac{\left(  N-k\right)  !\tau^{j-k}}{j\left(
N-j\right)  !}\nonumber\\
&  \hspace{0.5cm}+\frac{1}{\delta^{2}}\sum_{j=3}^{N}\sum_{k=3}^{j}\sum
_{s=1}^{k-3}\sum_{m=0}^{k-s-1}\frac{\left(  N-\left(  k-s-m\right)  \right)
!\left(  N-k\right)  !\tau^{j-k+m}}{j\left(  N-j\right)  !\left(  k-s\right)
\left(  N-\left(  k-s\right)  \right)  !} \label{f2_eps=0}%
\end{align}

\subsection{Asymptotics of $f_{0}$, $f_{1}$ and $f_{2}$ in $\varepsilon-$SIS
epidemics for large $N$}

\label{sec_coeff_eps_SIS_largeN}For large $N$, $\lim_{\varepsilon\rightarrow
0}f_{1}$ in (\ref{f1_eps=0}) behaves as%
\begin{align*}
\lim_{\varepsilon\rightarrow0}f_{1}  &  =F\left(  \tau\right)  =\frac
{1}{\delta}\sum_{j=1}^{N}\sum_{r=0}^{j-1}\frac{\left(  N-j+r\right)
!}{j\left(  N-j\right)  !}\tau^{r}\\
&  \sim\frac{1}{\delta}\sum_{j=1}^{N}\frac{1}{j}\sum_{r=0}^{j-1}N^{r}\tau
^{r}=\frac{1}{\delta}\sum_{j=1}^{N}\frac{1}{j}\frac{\left(  N\tau\right)
^{j}-1}{N\tau-1}%
\end{align*}
where we have used \cite[6.1.47]{Abramowitz}%
\begin{equation}
\frac{\Gamma\left(  N-a\right)  }{\Gamma\left(  N-b\right)  }=N^{a-b}\left(
1+\frac{\left(  a-b\right)  \left(  a+b-1\right)  }{2N}+O\left(
N^{-2}\right)  \right)  \label{Gamma_asymptotics}%
\end{equation}
to first order. The accurate asymptotic behavior of $F\left(  \tau\right)  $
is deduced in Appendix \ref{sec_order_F(tau_largeN}. Hence, using $x=N\tau$,
we find%
\begin{equation}
\delta\lim_{\varepsilon\rightarrow0}f_{1}=\frac{1+O\left(  \frac{1}{N}\right)
}{\left(  x-1\right)  }\sum_{j=1}^{N}\frac{x^{j}-1}{j} \label{f1_eps=0_largeN}%
\end{equation}
From (\ref{f1_eps=0_largeN}), the three regimes for $x$ lead to the following
growth. If $N\tau=x<1$, then%
\begin{align*}
\delta\lim_{\varepsilon\rightarrow0}f_{1}  &  =\frac{1+O\left(  \frac{1}%
{N}\right)  }{\left(  1-x\right)  }\left(  \sum_{j=1}^{N}\frac{1}{j}%
-\sum_{j=1}^{N}\frac{x^{j}}{j}\right) \\
&  =\frac{1+O\left(  \frac{1}{N}\right)  }{\left(  1-x\right)  }\left(
H_{N}+\log\left(  1-x\right)  \right)  =O\left(  \log N\right)
\end{align*}
If $x>1$, then%
\[
\delta\lim_{\varepsilon\rightarrow0}f_{1}=\frac{1+O\left(  \frac{1}{N}\right)
}{\left(  x-1\right)  }\sum_{j=1}^{N}\frac{x^{j}-1}{j}=O\left(  \frac{x^{N-1}%
}{N}\right)
\]
whereas for $x=1$%
\[
\delta\lim_{\varepsilon\rightarrow0}f_{1}=\left(  1+O\left(  \frac{1}%
{N}\right)  \right)  \sum_{j=1}^{N}1=N\left(  1+O\left(  \frac{1}{N}\right)
\right)  =O\left(  N\right)
\]

Invoking (\ref{Gamma_asymptotics}), the asymptotic expression of
$\lim_{\varepsilon\rightarrow0}f_{2}$ in (\ref{f1_eps=0}) for large $N$ is%
\[
\lim_{\varepsilon\rightarrow0}f_{2}\sim\frac{1}{2\delta^{2}}+\frac{1}%
{\delta^{2}}\sum_{j=3}^{N}\frac{\left(  N\tau\right)  ^{j-2}}{j}+\left(
\frac{\tau N+3}{2\delta^{2}}\right)  \sum_{j=3}^{N}\sum_{k=3}^{j}\frac{\left(
N\tau\right)  ^{j-k}}{j}+\frac{1}{\delta^{2}}\sum_{j=3}^{N}\sum_{k=3}^{j}%
\sum_{s=1}^{k-3}\sum_{m=0}^{k-s-1}\frac{\left(  N\tau\right)  ^{j-k+m}%
}{j\left(  k-s\right)  }%
\]
Using $x=N\tau$, we have%
\[
\delta^{2}\lim_{\varepsilon\rightarrow0}f_{2}\sim\frac{1}{2}+\left(
\frac{3x+1}{2\left(  x-1\right)  }\right)  \sum_{j=3}^{N}\frac{x^{j-2}}%
{j}-\frac{x+3}{2\left(  x-1\right)  }\sum_{j=3}^{N}\frac{1}{j}+R_{4}%
\]
where%
\[
R_{4}=\sum_{j=3}^{N}\sum_{k=3}^{j}\sum_{s=1}^{k-3}\sum_{m=0}^{k-s-1}%
\frac{x^{j-k+m}}{j\left(  k-s\right)  }%
\]
which can be simplified to%
\[
R_{4}=\frac{1}{\left(  x-1\right)  ^{2}}\sum_{j=3}^{N}\frac{x^{j}}{j}%
\sum_{m=3}^{j-1}\frac{1}{m}\left(  1-x^{m-j}-x^{-m}+x^{-j}\right)
\]
After further rearrangement of terms in $R_{4}$, we arrive at%
\begin{align*}
\delta^{2}\lim_{\varepsilon\rightarrow0}f_{2}  &  \sim\frac{1}{2}+\left(
\frac{3x+1}{2\left(  x-1\right)  }\right)  \sum_{j=3}^{N}\frac{x^{j-2}}%
{j}-\frac{x+3}{2\left(  x-1\right)  }\sum_{j=3}^{N}\frac{1}{j}\\
&  \hspace{0.5cm}+\frac{1}{\left(  x-1\right)  ^{2}}\sum_{j=3}^{N}\frac{x^{j}%
}{j}\left(  \sum_{m=3}^{j-1}\frac{1}{m}-\sum_{m=j+1}^{N}\frac{1}{m}\right)
+\frac{1}{2\left(  x-1\right)  ^{2}}\left(  \sum_{j=3}^{N}\frac{1}{j}\right)
^{2}\\
&  \hspace{0.5cm}-\frac{1}{\left(  x-1\right)  ^{2}}\sum_{j=1}^{N-3}%
\frac{x^{j}}{j}\left(  \sum_{m=3}^{3+j}\frac{1}{m}-\sum_{m=N+1-j}^{N}\frac
{1}{m}\right)
\end{align*}
If $x<1$, then%
\begin{align*}
\delta^{2}\lim_{\varepsilon\rightarrow0}f_{2}  &  \sim\frac{1}{2\left(
x-1\right)  ^{2}}\left(  \sum_{j=3}^{N}\frac{1}{j}\right)  ^{2}-\frac
{x+3}{2\left(  x-1\right)  }\sum_{j=3}^{N}\frac{1}{j}+O\left(  1\right) \\
&  =O\left(  \log^{2}N\right)
\end{align*}
whereas for $x>1$, collecting the largest power in $x$ yields
\[
\delta^{2}\lim_{\varepsilon\rightarrow0}f_{2}\sim\frac{x^{N-2}}{N}\left(
\frac{3}{2}+\sum_{m=3}^{N-1}\frac{1}{m}\right)  =O\left(  \frac{\log N}%
{N}x^{N-2}\right)
\]
For $x=1$, we can show that%
\[
\delta^{2}\lim_{\varepsilon\rightarrow0}f_{2}\sim\frac{N^{2}}{4}-\frac{9}%
{4}N+4\sum_{j=1}^{N}\frac{1}{j}-2=O\left(  N^{2}\right)
\]

\subsection{Scaling of $\zeta$ with $N$ in SIS epidemics (when $\varepsilon
\downarrow0$)}

\label{sec_zeta_SIS_three_regimes}In the limit for $\varepsilon\rightarrow0$,
the upper bound (\ref{upper_bound_zeta_first_order}) for $\zeta$ becomes, with
(\ref{f1_eps=0}),%
\begin{equation}
\zeta<-\frac{f_{0}}{f_{1}}=\frac{\delta}{\sum_{j=1}^{N}\frac{1}{j\left(
N-j\right)  !}\sum_{r=0}^{j-1}\left(  N-j+r\right)  !\tau^{r}}
\label{upper_bound_convergence_time_SIS}%
\end{equation}

Using the asymptotic expressions for $\lim_{\varepsilon\rightarrow0}f_{1}=1$,
$\lim_{\varepsilon\rightarrow0}f_{1}$ and $\lim_{\varepsilon\rightarrow0}%
f_{2}$, the second order Lagrange series
(\ref{zero_Lagrange_series_up_to_2_order}) for $\zeta$ is
\[
\zeta\approx-\frac{1}{f_{1}}-\frac{f_{2}}{f_{1}^{3}}%
\]
Thus, for $x=N\tau>1$,%
\[
-\zeta\approx O\left(  \frac{N}{x^{N-1}}+\frac{N^{2}\log N}{x^{2N-1}}\right)
=O\left(  \frac{N}{x^{N-1}}\right)
\]
illustrating that the first term in the Lagrange series is sufficient, leading
to our main result (\ref{zeta_SIS_KN_above_threshold_x>1}). Numerical
computations support this result. For $x<1$,%
\[
-\zeta\approx O\left(  \frac{1}{\log N}+\frac{\log^{2}N}{\log^{3}N}\right)
=O\left(  \frac{1}{\log N}\right)
\]
Since now the first and second term are of equal order, both need to be taken
into account. A second order Lagrange expansion is not sufficient and higher
order terms need to be evaluated in order to guarantee accuracy of $\zeta$. In
view of the dramatic increase in the computations, we refrain from pursuing
this track and content ourselves with numerical calculations. Finally, when
$x=1$, we have%
\[
-\zeta\approx O\left(  \frac{1}{N}+\frac{N^{2}}{N^{3}}\right)  =O\left(
\frac{1}{N}\right)
\]
leading to a similar conclusion as the case for $x<1$. In fact, we can compute
this zero a little more precise as%
\[
-\zeta\approx\left(  1+O\left(  \frac{1}{N}\right)  \right)  \left(
\frac{\delta}{N}\left(  1+\frac{\frac{N^{2}}{4}-\frac{9}{4}N+4\sum_{j=1}%
^{N}\frac{1}{j}-2}{N^{2}}\right)  \right)  =\frac{5\delta}{4N}\left(
1+O\left(  \frac{1}{N}\right)  \right)
\]
We also observe that $-\zeta$ is a rate, which is here naturally expressed in
units of the curing rate $\delta$.

\begin{figure}
[ptb]
\begin{center}
\includegraphics[
height=10.7525cm,
width=16.1298cm
]%
{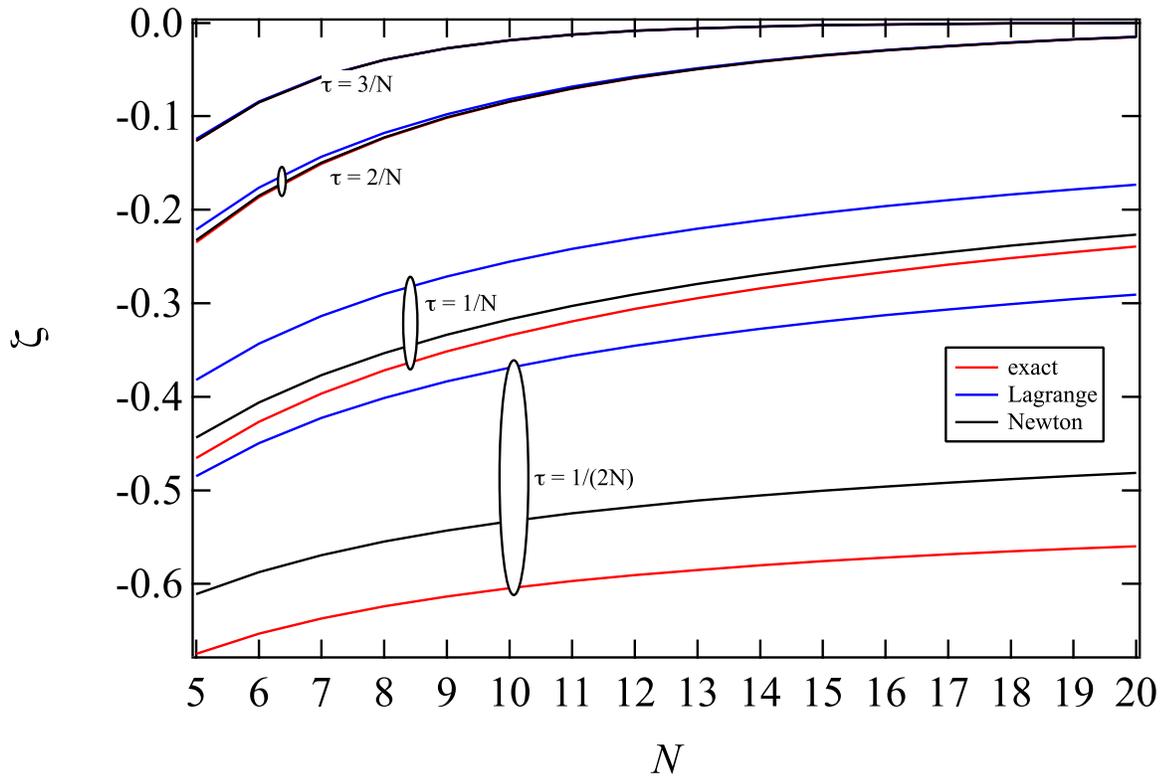}%
\caption{The second largest zero $\zeta$, exactly computed in red, by a second
order Lagrange series (\ref{zero_Lagrange_series_up_to_2_order}) in blue and
by Newton's identity (\ref{upper_bound_Newton_zeta}) in black as a function of
$N$ for $\varepsilon=10^{-5}$ and four values for the effective infection rate
$\tau=\left\{  1/(2N),1/N,2/N,3/N\right\}  $. }%
\label{Fig_zeta_eps1min5}%
\end{center}
\end{figure}

Fig.~\ref{Fig_zeta_eps1min5} shows the accuracy (for $\varepsilon=10^{-5}$) of
the second order Lagrange series (\ref{zero_Lagrange_series_up_to_2_order}),
the upper bound (\ref{upper_bound_Newton_zeta}) derived from the Newton
identities and the exact (numerical) computation of the second largest
eigenvalue of the infinitesimal generator matrix $Q$ of the continuous-time
$\varepsilon$-SIS Markov process on the complete graph $N$, for which the
epidemic threshold $\tau_{c}$ is slightly larger than $1/N$. These numerical
results confirm the order estimates (even for $\varepsilon\rightarrow0$)
above, at and below the epidemic threshold. Both Lagrange's second order and
Newton's upper bound are increasingly sharp for increasing values of $\tau$
above the epidemic threshold. Our exact asymptotics in
(\ref{zeta_SIS_KN_above_threshold_x>1}) of the order of $-\zeta=O\left(
\frac{N}{x^{N-1}}\right)  $ for $x>>1$ is difficult to verify for $N>10$ since
numerical root finders only provide an accuracy of about $10^{-10}$. For
$1<x<2.5$, (\ref{zeta_SIS_KN_above_threshold_x>1}) is verified up to $N=100$.
The relative accuracy for $\varepsilon<\frac{1}{N}$ is about the same as the
results shown in Fig.~\ref{Fig_zeta_eps1min5}.%

\section{Conclusion}

Our asymptotic order results agree with the general estimates of the average
lifetime of a SIS epidemics in Draief and Massouli\'{e}
\cite{Draief_Massoulie}. For large $t$, the probability of survival of the SIS
epidemics or probability that the life-time $T$ of an SIS epidemics exceeds
$t$ time units equals about%
\[
\Pr\left[  T>t\right]  \simeq e^{-\left\vert \zeta\right\vert t}%
\]
Hence, the life time of an epidemics (for large $t$) can be interpreted as
being exponentially distributed with mean $\frac{1}{\left\vert \zeta
\right\vert }$. In particular, above the epidemic threshold (equal to $x>1$ to
first order in $N$), the SIS epidemics in $K_{N}$ dies out exponentially in
time $t$ with decay rate $\zeta$, which tends to zero at least as fast as
$e^{-N\ln\frac{\tau}{\tau_{c}}}$, where $x=\frac{\tau}{\tau_{c}}>1$ is
measured in units of the epidemic threshold $\tau_{c}\sim\frac{1}{N}$ for
large $N$. This means that the probability that an SIS epidemic in any network
survives longer than $t$ time units is smaller than about $e^{-te^{-N\ln
\frac{\tau}{\tau_{c}}}}$ or that the average life time is at most $E\left[
T\right]  \simeq\frac{1}{\left\vert \zeta\right\vert }\simeq O\left(
e^{N\ln\frac{\tau}{\tau_{c}}}\right)  $, which is unrealistically long. Hence,
for sufficiently large $N$ and an effective infection rate $\tau>\tau_{c}$,
the SIS epidemics hardly ever dies in reality. When $\tau$ approaches
$\tau_{c}$, the decay rate $\zeta$ of the SIS epidemics decreases at least as
fast as $O\left(  \frac{1}{N}\right)  $, equivalent to an average life time
$E\left[  T\right]  \simeq O\left(  N\right)  $. Below the epidemic threshold
$x=\frac{\tau}{\tau_{c}}<1$, the decay rate $\zeta$ decreases at least as fast
as $O\left(  \frac{1}{\log N}\right)  $ and the average life time is about
$E\left[  T\right]  =O\left(  \log N\right)  $.

Finally, the lower bound results in \cite{Mountford2013} together with our
upper bound in (\ref{zeta_SIS_KN_above_threshold_x>1}) leads us to conclude
that for almost all graphs, the average time to absorption for $\tau>\tau_{c}$
is $E\left[  T\right]  =O\left(  e^{c_{G}N}\right)  $, where $c_{G}>0$. The
precise expression of $c_{G}$ for a given graph $G$ stays on the agenda of
future work.\medskip

\textbf{Acknowledgement} I am very grateful to Erik van Doorn for pointing me
to his and earlier work. Ruud van de Bovenkamp has provided me with numerical
data to test (\ref{zeta_SIS_KN_above_threshold_x>1}) for $N=100$ and various
$\tau>\tau_{c}$.

{\footnotesize
\bibliographystyle{plain}
\bibliography{cac,math,misc,net,pvm,qth,tel}
}

\appendix

\section{General tri-diagonal matrices}

\label{sec_general_tri-diagonal_matrices}We study the eigen-structure of
tri-diagonal matrices of the form%
\begin{equation}
P=\left[
\begin{array}
[c]{cccccccc}%
r_{0} & p_{0} & 0 & 0 & \cdots & 0 & 0 & 0\\
q_{1} & r_{1} & p_{1} & 0 & \cdots & 0 & 0 & 0\\
0 & q_{2} & r_{2} & p_{2} & \cdots & 0 & 0 & 0\\
\vdots & \vdots & \vdots & \vdots & \vdots & \vdots & \vdots & \vdots\\
0 & 0 & 0 & 0 & \cdots & q_{N-1} & r_{N-1} & p_{N-1}\\
0 & 0 & 0 & 0 & \cdots & 0 & q_{N} & r_{N}%
\end{array}
\right]  \label{P_general_tridiagonal_band_matrix}%
\end{equation}
where $p_{j}$ and $q_{j}$ are probabilities and where $P$ obeys the
stochasticity requirement $Pu=u$, where $u$ is the all-one vector. The matrix
$P$ frequently occurs in Markov theory, in particular, $P$ is the transition
probability matrix of the generalized random walk. The stochasticity
requirement reflects the fact that a Markov process must be in any of the
$N+1$ states. If $p_{j}=p$ and $q_{j}=q$, the matrix $P$ reduces to a Toeplitz
form for which the eigenvalues and eigenvectors can be explicitly written, as
shown in \cite{PVM_PerformanceAnalysisCUP}. Here, we consider the general
tri-diagonal matrix (\ref{P_general_tridiagonal_band_matrix}) and show how
orthogonal polynomials enter the scene. The theory for the discrete-time
generalized random walk is readily extended to that for the continuous-time
general birth-death process. While our approach is more algebraic, Karlin and
McGregor \cite{Karlin_McGregor1959} have presented a different, more
probabilistic and function-theoretic method, which is reviewed and
complemented by Van Doorn and Schrijner \cite{VanDoorn_Schrijner1995}.

When $P$ is written as a block matrix%
\[
P=\left[
\begin{array}
[c]{cc}%
A_{k\times k} & B_{k\times\left(  N+1-k\right)  }\\
C_{(N+1-k)\times k} & D_{\left(  N+1-k\right)  \times\left(  N+1-k\right)  }%
\end{array}
\right]
\]
then the matrix $B$ and $C$ only consist of one non-zero element, so that,
using the basic vector $e_{j}$ whose $j$-th component equals $1$ while all
others are zero, $B=p_{k-1}\left(  e_{k}\right)  _{k\times1}.\left(
e_{1}\right)  _{1\times\left(  N+1-k\right)  }^{T}$ and $C=q_{k}\left(
e_{1}\right)  _{\left(  N+1-k\right)  \times1}\left(  e_{k}\right)  _{1\times
k}^{T}$. The determinant of $P$ evaluated with Schur's formula
\[
\det\left[
\begin{array}
[c]{cc}%
A & B\\
C & D
\end{array}
\right]  =\det A\det\left(  D-CA^{-1}B\right)
\]
shows that
\begin{align*}
CA^{-1}B  &  =p_{k-1}q_{k}\left(  e_{1}\right)  _{\left(  N+1-k\right)
\times1}\left(  e_{k}\right)  _{1\times k}^{T}A^{-1}\left(  e_{k}\right)
_{k\times1}.\left(  e_{1}\right)  _{1\times\left(  N+1-k\right)  }^{T}\\
&  =p_{k-1}q_{k}\left(  A^{-1}\right)  _{kk}\left(  e_{1}\right)  _{\left(
N+1-k\right)  \times1}.\left(  e_{1}\right)  _{1\times\left(  N+1-k\right)
}^{T}%
\end{align*}
Thus, the matrix $CA^{-1}B$ only contains one non-zero element on position
$\left(  1,1\right)  $. Only the first element in $\widetilde{D}=D-CA^{-1}B$
is changed from $r_{k}$ in $D$ to $r_{k}-p_{k-1}q_{k}\left(  A^{-1}\right)
_{kk}$ in $\widetilde{D}$ and $p_{k-1}q_{k}\left(  A^{-1}\right)  _{kk}$ can
be considered as the coupling between the first $k-1$ states in the
generalized random walk and the remaining other states. If $p_{k-1}=0$ or/and
$q_{k}=0$, then $\det P=\det A\det D$, which is the product of two individual
tri-diagonal determinants. In that case, the Markov chain is reducible. Hence,
in the sequel, we assume that all elements of $P$ are non-zero and
time-independent so that the Markov chain is irreducible.

\subsection{A similarity transform}

We apply a similarity transform analogous to that of the Jacobi matrix for
orthogonal polynomials as studied in \cite[Section 10.6]{PVM_graphspectra}. If
there exists a similarity transform that makes the matrix $P$ symmetric, then
all eigenvalues of $P$ are real, because a similarity transform preserves the
eigenvalues. The simplest similarity transform is $H=$ diag$\left(
h_{1},h_{2},\ldots,h_{N+1}\right)  $ such that%
\[
\widetilde{P}=HPH^{-1}=\left[
\begin{array}
[c]{cccccccc}%
r_{0} & \frac{h_{1}}{h_{2}}p_{0} & 0 & 0 & \cdots & 0 & 0 & 0\\
\frac{h_{2}}{h_{1}}q_{1} & r_{1} & \frac{h_{2}}{h_{3}}p_{1} & 0 & \cdots & 0 &
0 & 0\\
0 & \frac{h_{3}}{h_{2}}q_{2} & r_{2} & \frac{h_{3}}{h_{4}}p_{2} & \cdots & 0 &
0 & 0\\
\vdots & \vdots & \vdots & \vdots & \vdots & \vdots & \vdots & \vdots\\
0 & 0 & 0 & 0 & \cdots & \frac{h_{N}}{h_{N-1}}q_{N-1} & r_{N-1} & \frac{h_{N}%
}{h_{N+1}}p_{N-1}\\
0 & 0 & 0 & 0 & \cdots & 0 & \frac{h_{N+1}}{h_{N}}q_{N} & r_{N}%
\end{array}
\right]
\]
Thus, in order to produce a symmetric matrix $\widetilde{P}=\widetilde{P}^{T}%
$, we need to require that $\left(  \widetilde{P}\right)  _{i,i-1}=\left(
\widetilde{P}\right)  _{i-1,i}$ for all $1\leq i\leq N$, implying that,%
\[
\frac{h_{i+1}}{h_{i}}q_{i}=\frac{h_{i}}{h_{i+1}}p_{i-1}%
\]
whence%
\[
\left(  \frac{h_{i+1}}{h_{i}}\right)  ^{2}=\frac{p_{i-1}}{q_{i}}%
\]
Assuming that all $p_{i}$ and $q_{i}$ are positive\footnote{If $p_{i-1}=0$ (or
$q_{i}=0$), then the states $0$ up to $i-1$ are uncoupled from the states $i$
up to $N$.}, we find that $h_{i+1}=\sqrt{\frac{p_{i-1}}{q_{i}}}h_{i}$ for
$1\leq i\leq N$ and we can choose $h_{1}=1$ such that
\begin{equation}
h_{i}=\sqrt{%
{\displaystyle\prod\limits_{k=1}^{i-1}}
\frac{p_{k-1}}{q_{k}}} \label{def_h_i}%
\end{equation}
and
\[
\frac{h_{i+1}}{h_{i}}q_{i}=\frac{h_{i}}{h_{i+1}}p_{i-1}=\sqrt{p_{i-1}q_{i}}%
\]
After the similarity transform $H$, the symmetric matrix $\widetilde{P}$
becomes%
\begin{equation}
\widetilde{P}=\left[
\begin{array}
[c]{cccccccc}%
r_{0} & \sqrt{p_{0}q_{1}} & 0 & 0 & \cdots & 0 & 0 & 0\\
\sqrt{p_{0}q_{1}} & r_{1} & \sqrt{p_{1}q_{2}} & 0 & \cdots & 0 & 0 & 0\\
0 & \sqrt{p_{1}q_{2}} & r_{2} & \sqrt{p_{2}q_{3}} & \cdots & 0 & 0 & 0\\
\vdots & \vdots & \vdots & \vdots & \vdots & \vdots & \vdots & \vdots\\
0 & 0 & 0 & 0 & \cdots & \sqrt{p_{N-2}q_{N-1}} & r_{N-1} & \sqrt{p_{N-1}q_{N}%
}\\
0 & 0 & 0 & 0 & \cdots & 0 & \sqrt{p_{N-1}q_{N}} & r_{N}%
\end{array}
\right]  \label{def_P_tilde}%
\end{equation}
In conclusion, if all $p_{i}$ and $q_{i}$ are positive, then all eigenvalues
of $P$ are real. Rather than solving the eigenvector $\widetilde{x}$ from the
eigenvalue equation $\widetilde{P}\widetilde{x}=\lambda\widetilde{x}$, we
determine the eigenvector $x$ as a function of $\lambda$ from the original
matrix $P$ for reasons explained below and use the similarity transform
$\widetilde{x}=Hx$, where $H$ is independent of $\lambda$, later for the
left-eigenvectors of $P$.

\subsection{Eigenvectors of $P$}

The right-eigenvector $x$ of $P$ belonging to eigenvalue $\lambda$ satisfies
$(P-\lambda I)x=0$ so that%
\[
\left\{
\begin{array}
[c]{rc}%
\left(  r_{0}-\lambda\right)  x_{0}+p_{0}x_{1}=0 & \\
q_{j}x_{j-1}+\left(  r_{j}-\lambda\right)  x_{j}+p_{j}x_{j+1}=0 & 1\leq j<N\\
q_{N}x_{N-1}+\left(  r_{N}-\lambda\right)  x_{N}=0 &
\end{array}
\right.
\]
We replace the last equation, that breaks the structure, by%
\[
q_{N}x_{N-1}+\left(  r_{N}-\lambda\right)  x_{N}+p_{N}x_{N+1}=0
\]
and the condition that $p_{N}x_{N+1}=0$. Using $r_{j}=1-q_{j}-p_{j}$ for
$0\leq j\leq N$ with $q_{0}=0$ and $p_{N}=0$ and making the dependence on
$\xi=\lambda-1$ explicit, the above set simplifies, subject to the condition
$p_{N}x_{N+1}\left(  \xi\right)  =0$, to%
\begin{equation}
\left\{
\begin{array}
[c]{lc}%
x_{1}\left(  \xi\right)  =\frac{p_{0}+\xi}{p_{0}}x_{0}\left(  \xi\right)  & \\
x_{j+1}\left(  \xi\right)  =\frac{p_{j}+q_{j}+\xi}{p_{j}}x_{j}\left(
\xi\right)  -\frac{q_{j}}{p_{j}}x_{j-1}\left(  \xi\right)  & 1\leq j<N
\end{array}
\right.  \label{difference_eigenvectors_P_grw}%
\end{equation}

For the stochastic matrix $P$, that obeys $Pu=u$, there holds that
$r_{N}=1-q_{N}$ so that $p_{N}=0$ and that the condition $p_{N}x_{N+1}=0$
seems to be obeyed. In the theory of orthogonal polynomials (see e.g.
\cite[Chapter 10]{PVM_graphspectra}), a similar trick is used where the
orthogonal polynomial $x_{N+1}\left(  \xi\right)  $ needs to vanish, because
$p_{N}$ is not necessarily zero in absence of the stochasticity requirement
$Pu=u$. The zeros of the orthogonal polynomial $x_{N+1}\left(  \xi\right)  $
are then equal to the eigenvalues of the corresponding Jabobi matrix.
Moreover, the powerful interlacing property for the zeros of the set $\left\{
x_{j}\left(  \xi\right)  \right\}  _{0\leq j\leq N+1}$ applies. We will return
to the condition $p_{N}x_{N+1}=0$ below.

Solving (\ref{difference_eigenvectors_P_grw}) iteratively for $\,j<N$,
\begin{align*}
x_{2}\left(  \xi\right)   &  =\frac{x_{0}\left(  \xi\right)  }{p_{0}p_{1}%
}\left(  \xi^{2}+\left(  q_{1}+p_{1}+p_{0}\right)  \xi+p_{1}p_{0}\right) \\
x_{3}\left(  \xi\right)   &  =\frac{x_{0}\left(  \xi\right)  }{p_{2}p_{1}%
p_{0}}\left(  \xi^{3}+\left(  q_{1}+q_{2}+p_{2}+p_{1}+p_{0}\right)  \xi
^{2}\right) \\
&  \hspace{0.5cm}+\left(  q_{2}q_{1}+q_{2}p_{0}+p_{2}q_{1}+p_{2}p_{1}%
+p_{2}p_{0}+p_{1}p_{0}\right)  \xi+p_{2}p_{1}p_{0}%
\end{align*}
reveals that $\frac{x_{j}\left(  \xi\right)  }{x_{0}\left(  \xi\right)  }$ is
a polynomial of degree $j$ in $\xi$ with positive coefficients, whose zeros
are all non-positive\footnote{If $P$ is a stochastic, irreducible matrix, then
the Perron-Frobenius Theorem \cite{PVM_graphspectra} states that the largest
(in absolute value) eigenvalue is one, hence $-1<$ $\lambda=\xi+1\leq1$.}.
This simple form is the main reason to consider the eigenvector components of
$P$ instead of $\widetilde{P}$. By inspection, the general form of $x_{j}%
(\xi)$ for $1\leq j\leq N$ is
\begin{equation}
x_{j}(\xi)=\frac{x_{0}\left(  \xi\right)  }{\prod_{m=0}^{j-1}p_{m}}\sum
_{k=0}^{j}c_{k}(j)\xi^{k} \label{general_eigenvector_component}%
\end{equation}
with\footnote{We use the convention that $\sum_{k=a}^{b}f\left(  k\right)  =0$
and $\prod\limits_{k=a}^{b}f\left(  k\right)  =1$ if $a>b$.}%
\begin{equation}%
\begin{array}
[c]{ccc}%
c_{j}(j)=1; & c_{j-1}(j)=\sum_{m=0}^{j-1}\left(  p_{m}+q_{m}\right)  ; &
c_{0}(j)=\prod_{m=0}^{j-1}p_{m};
\end{array}
\label{initial_values_c}%
\end{equation}
where $q_{0}=p_{N}=0$. By substituting (\ref{general_eigenvector_component})
into (\ref{difference_eigenvectors_P_grw}),
\[
\sum_{k=1}^{j-1}c_{k}(j+1)\xi^{k}=\sum_{k=1}^{j-1}\left[  \left(  q_{j}%
+p_{j}\right)  c_{k}(j)-q_{j}p_{j-1}c_{k}(j-1)+c_{k-1}(j)\right]  \xi^{k}%
\]
and equating the corresponding powers in $\xi$, a recursion relation for the
coefficients $c_{k}(j)$ for $0\leq k<j$ is obtained with $c_{j}\left(
j\right)  =1$,
\begin{equation}
c_{k}(j+1)=\left(  q_{j}+p_{j}\right)  c_{k}(j)-q_{j}p_{j-1}c_{k}%
(j-1)+c_{k-1}(j) \label{recursion_c_k(j)}%
\end{equation}
from which all coefficients can be determined as shown in Section
\ref{sec_solving_recursion_c_k(j)}. The stochasticity requirement $Pu=u$
implies that the right-eigenvector belonging to the largest eigenvalue
$\lambda=1$, equivalently to $\xi=\lambda-1=0$, equals $x\left(  0\right)  =u$.

We now express the left-eigenvectors of $P$ in terms of the right-eigenvector
by using the similarity transform $H$. Since $\widetilde{P}$ is symmetric, the
left- and right-eigenvectors are the same \cite[p. 222-223]{PVM_graphspectra}.
The left-eigenvector $x$ of $P$ equals $x=H^{-1}\widetilde{x}$, while the
right-eigenvector $y$ of $P$ equals $y=H\widetilde{x}$. Hence, we find that
$y=H^{2}x$ and explicitly with (\ref{general_eigenvector_component}) and
(\ref{def_h_i}),%
\begin{equation}
y_{j}(\xi)=\frac{y_{0}\left(  \xi\right)  }{\prod_{m=0}^{j-1}q_{m+1}}%
\sum_{k=0}^{j}c_{k}(j)\xi^{k} \label{general_left_eigenvector_component}%
\end{equation}

For any matrix, the left- and right-eigenvectors obey the orthogonality
equation%
\begin{equation}
x^{T}\left(  \xi\right)  y\left(  \xi^{\prime}\right)  =x^{T}\left(
\xi\right)  y\left(  \xi\right)  \delta_{\xi\xi^{\prime}}
\label{left-right-eigenvectors_orthogonality}%
\end{equation}
that holds for any pair of eigenvalues $\lambda=\xi+1$ and $\lambda^{\prime
}=\xi^{\prime}+1$ of that matrix. For symmetric matrices, usually, the
normalization
\begin{equation}
\widetilde{x}^{T}\left(  \xi\right)  \widetilde{x}\left(  \xi^{\prime}\right)
=\delta_{\xi\xi^{\prime}} \label{orthogonality_eigenvectors_P}%
\end{equation}
is chosen, which implies, after the similarity transform $H=$ diag$\left(
h_{i}\right)  $, that $x^{T}\left(  \xi\right)  y\left(  \xi^{\prime}\right)
=\widetilde{x}^{T}\left(  \xi\right)  \widetilde{x}\left(  \xi^{\prime
}\right)  =\delta_{\xi\xi^{\prime}}$ and that $x^{T}\left(  \xi\right)
H^{2}x\left(  \xi^{\prime}\right)  =\delta_{\xi\xi^{\prime}}$ and similarly
that $y^{T}\left(  \xi\right)  H^{-2}y\left(  \xi^{\prime}\right)
=\delta_{\xi\xi^{\prime}}$. These normalizations of the eigenvector components
imply, using (\ref{general_eigenvector_component}) and
(\ref{general_left_eigenvector_component}) and with the definition of the
polynomial for $0\leq j\leq N+1$%
\begin{equation}
\rho_{j}\left(  \xi\right)  =\sum_{k=0}^{j}c_{k}(j)\xi^{k}
\label{def_rj_orthogonal_polynomial_of_P}%
\end{equation}
that
\begin{equation}
x_{0}\left(  \xi\right)  y_{0}\left(  \xi\right)  =x_{0}^{2}\left(
\xi\right)  =y_{0}^{2}\left(  \xi\right)  =\left(  1+\sum_{j=1}^{N}\frac
{\rho_{j}^{2}\left(  \xi\right)  }{\prod_{m=0}^{j-1}q_{m+1}p_{m}}\right)
^{-1} \label{normalization_x0_y_0}%
\end{equation}
The orthogonality equation (\ref{left-right-eigenvectors_orthogonality})
together with our choice of normalization, $\widetilde{x}^{T}\left(
\xi\right)  \widetilde{x}\left(  \xi^{\prime}\right)  =\delta_{\xi\xi^{\prime
}}$, lead to a couple of important consequences.

Using the general form (\ref{general_eigenvector_component}), the initially
made condition $p_{N}x_{N+1}=0$ translates to%
\[
p_{N}x_{N+1}(\xi)=\frac{x_{0}\left(  \xi\right)  \rho_{N+1}\left(  \xi\right)
}{\prod_{m=0}^{N-1}p_{m}}=0
\]
Since $\rho_{j}\left(  \xi\right)  $ is a polynomial,
(\ref{normalization_x0_y_0}) indicates that neither $x_{0}\left(  \xi\right)
$ nor $y_{0}\left(  \xi\right)  $ can vanish for finite $\xi$ so that the
initial condition is met provided%
\begin{equation}
\rho_{N+1}\left(  \xi\right)  =0 \label{zeros_x_N+1=eigenvalues_P}%
\end{equation}
which closely corresponds to results in the theory of orthogonal polynomials.
Thus, $\rho_{j}\left(  \xi\right)  $ should be considered as orthogonal
polynomial, rather than $x_{j}\left(  \xi\right)  $ due to the scaling of
$x_{0}\left(  \xi\right)  $, defined in (\ref{normalization_x0_y_0}). For the
set of orthogonal polynomials $\left\{  \rho_{j}\left(  \xi\right)  \right\}
_{0\leq j\leq N+1}$ \emph{interlacing} applies, which means that the zeros of
$\rho_{j}\left(  \xi\right)  $ interlace with those of $\rho_{l}\left(
\xi\right)  $ for all $1\leq l\neq j\leq N+1$. Moreover, the eigenvalues of
$P$ are equal to the zeros of $\rho_{N+1}\left(  \xi\right)  $ in
(\ref{zeros_x_N+1=eigenvalues_P}).

For stochastic matrices, the left-eigenvector $y(0)$ belonging to $\xi=0$
equals the steady-state vector $\pi$ (see \cite{PVM_PerformanceAnalysisCUP}).
For $\xi=0$, the orthogonality relation
(\ref{left-right-eigenvectors_orthogonality}) becomes $u^{T}y\left(
\xi^{\prime}\right)  =0$ and $u^{T}y\left(  0\right)  =u^{T}\pi=1$, from which
the $j$-th component in (\ref{general_left_eigenvector_component}) of the
left-eigenvector $y\left(  0\right)  =\pi$ follows, for $1\leq j\leq N$, as%
\begin{equation}
\pi_{j}=y_{j}(0)=\frac{y_{0}\left(  0\right)  c_{0}\left(  j\right)  }%
{\prod_{m=0}^{j-1}q_{m+1}}=\frac{\prod_{m=0}^{j-1}\frac{p_{m}}{q_{m+1}}%
}{1+\sum_{k=1}^{N}\prod_{m=0}^{k-1}\frac{p_{m}}{q_{m+1}}}
\label{steady_state_general_random_walk}%
\end{equation}
which precisely equal the well-known steady-state probabilities of the
generalized random walk \cite[p. 207]{PVM_PerformanceAnalysisCUP}. For
$\xi\neq0$, the orthogonality relation
(\ref{left-right-eigenvectors_orthogonality}) and the fact that $y_{0}\left(
\xi\right)  $ is non-zero for finite $\xi$ imply that%

\[
0=1+\sum_{j=1}^{N}\frac{1}{\prod_{m=0}^{j-1}q_{m+1}}\sum_{k=0}^{j}c_{k}%
(j)\xi^{k}=1+\sum_{j=1}^{N}\prod_{m=0}^{j-1}\frac{p_{m}}{q_{m+1}}+\sum
_{k=1}^{N}\left(  \sum_{j=k}^{N}\frac{c_{k}(j)}{\prod_{m=0}^{j-1}q_{m+1}%
}\right)  \xi^{k}%
\]
We write the right-hand side polynomial as%
\begin{equation}
\sum_{k=0}^{N}f_{k}\xi^{k}=f_{N}\prod\limits_{k=1}^{N}\left(  \xi
-z_{k}\right)  \label{characteristic_polynomial_P_grw_in_zeros}%
\end{equation}
where $f_{0}=\frac{1}{\pi_{0}}$ by (\ref{steady_state_general_random_walk})
and where, for $k>0$,%
\begin{equation}
f_{k}=\sum_{j=k}^{N}\frac{c_{k}(j)}{\prod_{m=0}^{j-1}q_{m+1}} \label{def_f_k}%
\end{equation}
and, explicitly,%
\begin{align*}
f_{N}  &  =\frac{1}{\prod_{m=0}^{N-1}q_{m+1}}\\
f_{N-1}  &  =\frac{1}{\prod_{m=0}^{N-2}q_{m+1}}+\frac{\sum_{m=0}^{N-1}\left(
p_{m}+q_{m}\right)  }{\prod_{m=0}^{N-1}q_{m+1}}%
\end{align*}
Relation (\ref{def_f_k}) illustrates that all coefficients $f_{k}$ are
non-negative. Moreover, the orthogonality relation
(\ref{left-right-eigenvectors_orthogonality}) implies that the polynomial
$\sum_{k=0}^{N}f_{k}\xi^{k}$ possesses the same zeros as $\frac{c_{P}\left(
\xi\right)  }{\xi}$, where
\[
c_{P}\left(  \xi\right)  =\det\left(  P-\left(  \xi+1\right)  I\right)
=\xi\prod\limits_{k=1}^{N}\left(  z_{k}-\xi\right)
\]
is the $N+1$ degree characteristic polynomial of the matrix $P$, in
particular,%
\begin{equation}
\frac{c_{P}\left(  \xi\right)  }{\left(  -1\right)  ^{N}\xi}=\frac{1}{f_{N}%
}\sum_{k=0}^{N}f_{k}\xi^{k} \label{characteristic_pol_P}%
\end{equation}
Finally, since the eigenvalues of $P$ also obey
(\ref{zeros_x_N+1=eigenvalues_P}) so that
\[
\frac{c_{P}\left(  \xi\right)  }{\left(  -1\right)  ^{N}\xi}=\frac{1}{f_{N}%
}\sum_{k=0}^{N}f_{k}\xi^{k}=\sum_{k=1}^{N+1}c_{k}(N+1)\xi^{k-1}%
\]
after equating corresponding powers in $\xi$, we find that%
\begin{equation}
c_{k+1}(N+1)=\frac{f_{k}}{f_{N}}=\sum_{j=k}^{N}c_{k}(j)\prod_{m=j}%
^{N-1}q_{m+1} \label{coefficients_c_k(N+1)}%
\end{equation}
In summary, the stochasticity property of $P$ provides us with an additional
relation (\ref{coefficients_c_k(N+1)}) on the coefficients of $c_{P}\left(
\xi\right)  $, that is not necessarily obeyed for general orthogonal polynomials.

\subsection{Solving the recursion (\ref{recursion_c_k(j)})}

\label{sec_solving_recursion_c_k(j)}We now propose two different types of
solutions of the recursion (\ref{recursion_c_k(j)}) for the coefficients
$c_{k}\left(  j\right)  $ of $x_{j}\left(  \xi\right)  $ in
(\ref{general_eigenvector_component}).

\begin{theorem}
\label{theorem_solution_c_k(j)_for_k=j-m} A recursion relation for
$c_{j-m}\left(  j\right)  $, valid for $2\leq m\leq j$, is%
\begin{equation}
c_{j-m}(j)=\sum_{l=0}^{j-1}\left(  \left(  q_{l}+p_{l}\right)  c_{l-m+1}%
(l)-q_{l}p_{l-1}c_{l-m+1}(l-1)\right)  \label{t_m(j)}%
\end{equation}

\end{theorem}

\textbf{Proof}: Letting $k=j-m$ in (\ref{recursion_c_k(j)}) yields%
\[
c_{j+1-(m+1)}(j+1)=\left(  q_{j}+p_{j}\right)  c_{j-m}(j)-q_{j}p_{j-1}%
c_{j-1-(m-1)}(j-1)+c_{j-(m+1)}(j)
\]
With $t_{m}\left(  j\right)  =c_{j-m}(j)$, the above equation transforms into
the difference equation%
\[
t_{m+1}\left(  j+1\right)  =t_{m+1}\left(  j\right)  +\left(  q_{j}%
+p_{j}\right)  t_{m}\left(  j\right)  -q_{j}p_{j-1}t_{m-1}\left(  j-1\right)
\]
whose solution is%
\[
t_{m+1}\left(  j\right)  =\sum_{l=0}^{j-1}\left(  \left(  q_{l}+p_{l}\right)
t_{m}\left(  l\right)  -q_{l}p_{l-1}t_{m-1}\left(  l-1\right)  \right)
\]
With the initial values $t_{0}\left(  j\right)  =1$ and $t_{1}\left(
j\right)  =\sum_{m=0}^{j-1}\left(  p_{m}+q_{m}\right)  $ from
(\ref{initial_values_c}), all $t_{m}\left(  j\right)  $ can be iteratively
found from (\ref{t_m(j)}).\hfill$\square\medskip$

Thus, letting $m=2$ in (\ref{t_m(j)}) yields%
\begin{equation}
c_{j-2}(j)=\sum_{l=0}^{j-1}\left(  \left(  q_{l}+p_{l}\right)  \sum
_{m=0}^{l-1}\left(  p_{m}+q_{m}\right)  -q_{l}p_{l-1}\right)
\label{c_k(j)_for_k=j-2}%
\end{equation}
Next, for $m=3$ in (\ref{t_m(j)}), we have%
\[
c_{j-3}(j)=\sum_{l=0}^{j-1}\left(  \left(  q_{l}+p_{l}\right)  \sum_{l_{1}%
=0}^{l-1}\left(  \left(  q_{l_{1}}+p_{l_{1}}\right)  \sum_{m=0}^{l_{1}%
-1}\left(  p_{m}+q_{m}\right)  -q_{l_{1}}p_{l_{1}-1}\right)  -q_{l}p_{l-1}%
\sum_{m=0}^{l-2}\left(  p_{m}+q_{m}\right)  \right)
\]
and so on.

For the polynomial in (\ref{characteristic_polynomial_P_grw_in_zeros}), the
next general expression will prove more useful.

\begin{theorem}
\label{theorem_solution_recursion_c_k(j)} The explicit general expression for
the coefficients $c_{k}\left(  j\right)  $ in terms of $c_{k-1}\left(
l\right)  $ for all $l\geq k-1$ is%
\begin{align}
c_{k}\left(  j\right)   &  =\prod\limits_{m=k}^{j-1}p_{m}+\sum_{l=0}%
^{j-k-1}\prod\limits_{m=k}^{j-l-1}q_{m}\prod\limits_{m=j-l}^{j-1}p_{m}%
+\sum_{m=0}^{k-1}\left(  p_{m}+q_{m}\right)  \sum_{l=0}^{j-k-1}\prod
\limits_{m=k+1}^{j-l-1}q_{m}\prod\limits_{m=j-l}^{j-1}p_{m}\nonumber\\
&  \hspace{0.5cm}+\sum_{l=0}^{j-k-1}\sum_{s=1}^{j-l-k-1}c_{k-1}\left(
j-l-s\right)  \prod\limits_{m=j-l+1-s}^{j-l-1}q_{m}\prod\limits_{m=j-l}%
^{j-1}p_{m} \label{ck(j)_as_sum_c_(k-1)(x)}%
\end{align}

\end{theorem}

\textbf{Proof:} Rewriting (\ref{recursion_c_k(j)}) as%
\[
c_{k}(j+1)-p_{j}c_{k}(j)=q_{j}\left\{  c_{k}(j)-p_{j-1}c_{k}(j-1)\right\}
+c_{k-1}(j)
\]
and defining $b_{k}\left(  j\right)  =c_{k}(j)-p_{j-1}c_{k}(j-1)$ shows that
the second order recursion (\ref{recursion_c_k(j)}) in $j$ can be decomposed
into two first order recursions in $j$%
\[
\left\{
\begin{array}
[c]{l}%
c_{k}(j)=p_{j-1}c_{k}(j-1)+b_{k}\left(  j\right) \\
b_{k}\left(  j\right)  =q_{j-1}b_{k}\left(  j-1\right)  +c_{k-1}(j-1)
\end{array}
\right.
\]
Since $k<j$, the choice for $j=k+1$ yields%
\begin{align*}
b_{k}\left(  k+1\right)   &  =c_{k}(k+1)-p_{k}c_{k}(k)\\
&  =\sum_{m=0}^{k}\left(  p_{m}+q_{m}\right)  -p_{k}=q_{k}+\sum_{m=0}%
^{k-1}\left(  p_{m}+q_{m}\right)
\end{align*}
Iterating the first recursion downwards yields%
\begin{align*}
c_{k}(j)  &  =p_{j-1}p_{j-2}c_{k}(j-2)+p_{j-1}b_{k}\left(  j-1\right)
+b_{k}\left(  j\right) \\
&  =p_{j-1}p_{j-2}p_{j-3}c_{k}(j-3)+p_{j-1}p_{j-2}b_{k}\left(  j-2\right)
+p_{j-1}b_{k}\left(  j-1\right)  +b_{k}\left(  j\right) \\
&  =p_{j-1}p_{j-2}p_{j-3}p_{j-4}c_{k}(j-4)+p_{j-1}p_{j-2}p_{j-3}b_{k}\left(
j-3\right)  +p_{j-1}p_{j-2}b_{k}\left(  j-2\right)  +p_{j-1}b_{k}\left(
j-1\right)  +b_{k}\left(  j\right)
\end{align*}
from which we deduce that%
\[
c_{k}\left(  j\right)  =c_{k}(j-p)\prod\limits_{m=j-p}^{j-1}p_{m}+\sum
_{l=0}^{p-1}b_{k}\left(  j-l\right)  \prod\limits_{m=j-l}^{j-1}p_{m}%
\]
When $j-p=k$, then $c_{k}\left(  k\right)  =1$ and thus%
\begin{equation}
c_{k}\left(  j\right)  =\prod\limits_{m=k}^{j-1}p_{m}+\sum_{l=0}^{j-k-1}%
b_{k}\left(  j-l\right)  \prod\limits_{m=j-l}^{j-1}p_{m}
\label{ck(j)_in_sum_bk}%
\end{equation}
Similarly, we iterate the second recursion downwards,%
\begin{align*}
b_{k}\left(  j\right)   &  =q_{j-1}q_{j-2}b_{k}\left(  j-2\right)
+q_{j-1}c_{k-1}(j-2)+c_{k-1}(j-1)\\
&  =q_{j-1}q_{j-2}q_{j-3}b_{k}\left(  j-3\right)  +q_{j-1}q_{j-2}%
c_{k-1}(j-3)+q_{j-1}c_{k-1}(j-2)+c_{k-1}(j-1)
\end{align*}
which suggests that%
\[
b_{k}\left(  j\right)  =b_{k}(j-p)\prod\limits_{m=j-p}^{j-1}q_{m}+\sum
_{l=1}^{p}c_{k-1}\left(  j-l\right)  \prod\limits_{m=j+1-l}^{j-1}q_{m}%
\]
For $j-p=k+1$ or $p=j-k-1$, we have%
\begin{align}
b_{k}\left(  j\right)   &  =b_{k}(k+1)\prod\limits_{m=k+1}^{j-1}q_{m}%
+\sum_{l=1}^{j-k-1}c_{k-1}\left(  j-l\right)  \prod\limits_{m=j+1-l}%
^{j-1}q_{m}\nonumber\\
&  =\prod\limits_{m=k}^{j-1}q_{m}+\prod\limits_{m=k+1}^{j-1}q_{m}\sum
_{m=0}^{k-1}\left(  p_{m}+q_{m}\right)  +\sum_{l=1}^{j-k-1}c_{k-1}\left(
j-l\right)  \prod\limits_{m=j+1-l}^{j-1}q_{m} \label{bk(j)_in_c_k-1}%
\end{align}
Combining (\ref{ck(j)_in_sum_bk}) and (\ref{bk(j)_in_c_k-1}) yields
(\ref{ck(j)_as_sum_c_(k-1)(x)}).\hfill$\square\medskip$

For $k=1$ and using $c_{0}(j)=\prod_{m=0}^{j-1}p_{m}$, we find from
(\ref{ck(j)_as_sum_c_(k-1)(x)}) that%
\begin{equation}
c_{1}\left(  j\right)  =\sum_{l=0}^{j-1}\sum_{s=0}^{j-1-l}\prod\limits_{m=0}%
^{j-2-l-s}p_{m}\prod\limits_{m=j-l-s}^{j-1-l}q_{m}\prod\limits_{m=j-l}%
^{j-1}p_{m} \label{c1(j)_explicit}%
\end{equation}
Introducing the expression (\ref{c1(j)_explicit}) for $c_{1}\left(  j\right)
$ into (\ref{ck(j)_as_sum_c_(k-1)(x)}) produces the explicit form for
$c_{2}\left(  j\right)  $,%
\begin{align}
c_{2}\left(  j\right)   &  =\prod\limits_{m=2}^{j-1}p_{m}+\left(  p_{0}%
+p_{1}+q_{1}+q_{2}\right)  \sum_{l=0}^{j-3}\prod\limits_{m=3}^{j-l-1}%
q_{m}\prod\limits_{m=j-l}^{j-1}p_{m}\nonumber\\
&  \hspace{0.5cm}+\sum_{l=0}^{j-3}\sum_{s=1}^{j-l-3}\sum_{l_{1}=0}%
^{j-l-s-1}\sum_{l_{2}=0}^{j-l-s-l_{1}-1}\prod\limits_{m=0}^{j-l-s-2-l_{1}%
-l_{2}}p_{m}\prod\limits_{m=j-l-s-l_{1}-l_{2}}^{j-l-s-1-l_{1}}q_{m}%
\prod\limits_{m=j-l-s-l_{1}}^{j-l-s-1}p_{m}\prod\limits_{m=j-l+1-s}%
^{j-l-1}q_{m}\prod\limits_{m=j-l}^{j-1}p_{m} \label{c2(j)_explicit}%
\end{align}
and so on. In this way, all coefficients $c_{k}\left(  j\right)  $ in the
polynomial (\ref{general_eigenvector_component}) can be explicitly
determined\footnote{For $j=k$ in (\ref{ck(j)_as_sum_c_(k-1)(x)}), we find
indeed that $c_{k}\left(  k\right)  =1$ (based on our convention).}. Since all
$p_{j}$ and $q_{j}$ are probabilities and thus non-negative, the recursion
(\ref{ck(j)_as_sum_c_(k-1)(x)}) together with $c_{0}(j)=\prod_{m=0}^{j-1}%
p_{m}$ illustrates that all coefficients $c_{k}\left(  j\right)  $ are non-negative.

\subsection{The second set of orthogonality conditions}

Since the matrix $\widetilde{X}$, with the eigenvectors $\widetilde{x}$ of the
symmetric matrix $\widetilde{P}$ as columns, is orthogonal, it holds that%
\[
\widetilde{X}^{T}\widetilde{X}=\widetilde{X}\widetilde{X}^{T}=I
\]
and the last equation means that%
\[
\sum_{\lambda\in\left\{  \lambda_{1},\lambda_{2},\ldots,\lambda_{N}%
,\lambda_{N+1}\right\}  }\widetilde{x}_{j}\left(  \lambda\right)
\widetilde{x}_{m}\left(  \lambda\right)  =\delta_{jm}%
\]
where $\lambda_{1}=1\geq\lambda_{2}\geq\cdots\geq\lambda_{N+1}$ are the
eigenvalues of $P$ corresponding to the zeros of $c_{P}\left(  \xi\right)  $
by $\xi=\lambda-1$ and%
\[
\sum_{\lambda\in\left\{  \lambda_{1},\lambda_{2},\ldots,\lambda_{N}%
,\lambda_{N+1}\right\}  }h_{j+1}h_{m+1}x_{j}\left(  \lambda\right)
x_{m}\left(  \lambda\right)  =\delta_{jm}%
\]
Using (\ref{general_eigenvector_component}) and (\ref{def_h_i}) yields%
\[
\sum_{\lambda\in\left\{  \lambda_{1},\lambda_{2},\ldots,\lambda_{N}%
,\lambda_{N+1}\right\}  }x_{0}^{2}\left(  \lambda-1\right)  \sum_{k=0}%
^{j}c_{k}(j)\left(  \lambda-1\right)  ^{k}\sum_{l=0}^{k}c_{l}(k)\left(
\lambda-1\right)  ^{l}=\sqrt{%
{\displaystyle\prod\limits_{m=1}^{j}}
q_{m}%
{\displaystyle\prod\limits_{m=1}^{k}}
q_{m}}\delta_{jk}%
\]
which we rewrite, with the definition (\ref{def_rj_orthogonal_polynomial_of_P}%
), as%
\[
\sum_{\lambda\in\left\{  \lambda_{1},\lambda_{2},\ldots,\lambda_{N}%
,\lambda_{N+1}\right\}  }x_{0}^{2}\left(  \lambda-1\right)  \rho_{j}\left(
\lambda-1\right)  \rho_{k}\left(  \lambda-1\right)  =\delta_{jk}%
{\displaystyle\prod\limits_{m=1}^{k}}
q_{m}%
\]
Finally, introducing the Dirac delta-function, the left-hand side is rewritten
as an integral%
\begin{align*}
I  &  =\sum_{\lambda\in\left\{  \lambda_{1},\lambda_{2},\ldots,\lambda
_{N},\lambda_{N+1}\right\}  }x_{0}^{2}\left(  \lambda-1\right)  \rho
_{j}\left(  \lambda-1\right)  \rho_{k}\left(  \lambda-1\right) \\
&  =\sum_{j=1}^{N+1}\int_{-1}^{1}\delta\left(  \lambda-\lambda_{j}\right)
x_{0}^{2}\left(  \lambda-1\right)  \rho_{j}\left(  \lambda-1\right)  \rho
_{k}\left(  \lambda-1\right)  d\lambda
\end{align*}
because the eigenvalues of $P$ lie between $\left[  -1,1\right]  $. Further,
\begin{align*}
I  &  =\int_{-1}^{1}d\lambda x_{0}^{2}\left(  \lambda-1\right)  \rho
_{j}\left(  \lambda-1\right)  \rho_{k}\left(  \lambda-1\right)  \delta\left(
\det\left(  P-\lambda I\right)  \right)  \left\vert \left.  \frac{d\det\left(
P-xI\right)  }{dx}\right\vert _{x=\lambda}\right\vert \\
&  =\int_{-2}^{0}d\xi x_{0}^{2}\left(  \xi\right)  \rho_{j}\left(  \xi\right)
\rho_{k}\left(  \xi\right)  \delta\left(  \det\left(  P-\left(  \xi+1\right)
I\right)  \right)  \left\vert \left.  \frac{d\det\left(  P-xI\right)  }%
{dx}\right\vert _{x=\xi+1}\right\vert
\end{align*}
Defining the weight function as%
\begin{align*}
w\left(  \xi\right)   &  =x_{0}^{2}\left(  \xi\right)  \delta\left(
\det\left(  P-\left(  \xi+1\right)  I\right)  \right)  \left\vert \left.
\frac{d\det\left(  P-xI\right)  }{dx}\right\vert _{x=\xi+1}\right\vert \\
&  =x_{0}^{2}\left(  \xi\right)  \delta\left(  c_{P}\left(  \xi\right)
\right)  \left\vert \frac{dc_{P}\left(  \xi\right)  }{d\xi}\right\vert
=\sum_{j=1}^{N+1}x_{0}^{2}\left(  \xi_{j}\right)  \delta\left(  \xi-\xi
_{j}\right)
\end{align*}
we finally obtain the orthogonality condition for the orthogonal polynomials
$r_{j}$ and $r_{k}$ as%
\[
\int_{-2}^{0}w\left(  \xi\right)  \rho_{j}\left(  \xi\right)  \rho_{k}\left(
\xi\right)  d\xi=\delta_{jk}%
{\displaystyle\prod\limits_{m=1}^{k}}
q_{m}%
\]
In summary, the derivation provides an explicit way to determine the weight
function $w\left(  \xi\right)  $ in the orthogonality relation corresponding
to a tri-diagonal stochastic matrix $P$.

\subsection{The Christoffel-Darboux formula for eigenvectors of $P$}

We derive the Christoffel-Darboux formula (see \cite[p. 357]{PVM_graphspectra}%
) for the matrix $P$. Indeed, multiply the equation for $x_{j+1}\left(
\xi\right)  $ in (\ref{difference_eigenvectors_P_grw}) by $x_{j}\left(
\omega\right)  $%
\[
p_{j}x_{j+1}\left(  \xi\right)  x_{j}\left(  \omega\right)  =\xi x_{j}\left(
\xi\right)  x_{j}\left(  \omega\right)  +\left(  p_{j}+q_{j}\right)
x_{j}\left(  \xi\right)  x_{j}\left(  \omega\right)  -q_{j}x_{j-1}\left(
\xi\right)  x_{j}\left(  \omega\right)
\]
Letting $\xi\rightarrow\omega$ in (\ref{difference_eigenvectors_P_grw}) and
multiply both sides by $x_{j}\left(  \xi\right)  $,
\[
p_{j}x_{j+1}\left(  \omega\right)  x_{j}\left(  \xi\right)  =\omega
x_{j}\left(  \xi\right)  x_{j}\left(  \omega\right)  +\left(  p_{j}%
+q_{j}\right)  x_{j}\left(  \xi\right)  x_{j}\left(  \omega\right)
-q_{j}x_{j}\left(  \xi\right)  x_{j-1}\left(  \omega\right)
\]
Subtracting both equation yields,%
\[
p_{j}\left\{  x_{j+1}\left(  \xi\right)  x_{j}\left(  \omega\right)
-x_{j+1}\left(  \omega\right)  x_{j}\left(  \xi\right)  \right\}
+q_{j}\left\{  x_{j-1}\left(  \xi\right)  x_{j}\left(  \omega\right)
-x_{j}\left(  \xi\right)  x_{j-1}\left(  \omega\right)  \right\}  =\left(
\xi-\omega\right)  x_{j}\left(  \xi\right)  x_{j}\left(  \omega\right)
\]
Now, we transform to $x_{j}\left(  \xi\right)  =\frac{\widetilde{x}_{j}\left(
\xi\right)  }{h_{j+1}}$,%
\[
\frac{p_{j}}{h_{j+2}h_{j+1}}\left\{  \widetilde{x}_{j+1}\left(  \xi\right)
\widetilde{x}_{j}\left(  \omega\right)  -\widetilde{x}_{j+1}\left(
\omega\right)  \widetilde{x}_{j}\left(  \xi\right)  \right\}  +\frac{q_{j}%
}{h_{j}h_{j+1}}\left\{  \widetilde{x}_{j-1}\left(  \xi\right)  \widetilde
{x}_{j}\left(  \omega\right)  -\widetilde{x}_{j-1}\left(  \omega\right)
\widetilde{x}_{j}\left(  \xi\right)  \right\}  =\frac{\left(  \xi
-\omega\right)  }{h_{j+1}^{2}}\widetilde{x}_{j}\left(  \xi\right)
\widetilde{x}_{j}\left(  \omega\right)
\]
Using (\ref{def_h_i}) shows that $\frac{p_{j}}{h_{j+2}h_{j+1}}=\frac
{\sqrt{p_{j}q_{j+1}}}{h_{j+1}^{2}}$ and $\frac{q_{j}}{h_{j}h_{j+1}}%
=\frac{\sqrt{p_{j-1}q_{j}}}{h_{j+1}^{2}}$ so that%
\[
g_{j+1}-g_{j}=\left(  \xi-\omega\right)  \widetilde{x}_{j}\left(  \xi\right)
\widetilde{x}_{j}\left(  \omega\right)
\]
where%
\[
g_{j}=\sqrt{p_{j-1}q_{j}}\left\{  \widetilde{x}_{j-1}\left(  \omega\right)
\widetilde{x}_{j}\left(  \xi\right)  -\widetilde{x}_{j-1}\left(  \xi\right)
\widetilde{x}_{j}\left(  \omega\right)  \right\}
\]
Summing over $j\in\left[  0,m\right]  $,%
\[
\left(  \xi-\omega\right)  \sum_{j=0}^{m}\widetilde{x}_{j}\left(  \xi\right)
\widetilde{x}_{j}\left(  \omega\right)  =\sum_{j=0}^{m}g_{j+1}-\sum_{j=0}%
^{m}g_{j}=g_{m+1}-g_{0}%
\]
where $g_{0}=0$ because $\widetilde{x}_{-1}=0$. Hence, we arrive at the
Christoffel-Darboux sum for the eigenvectors of $\widetilde{P}$,%
\[
\left(  \xi-\omega\right)  \sum_{j=0}^{m}\widetilde{x}_{j}\left(  \xi\right)
\widetilde{x}_{j}\left(  \omega\right)  =\sqrt{p_{m}q_{m+1}}\left\{
\widetilde{x}_{m}\left(  \omega\right)  \widetilde{x}_{m+1}\left(  \xi\right)
-\widetilde{x}_{m}\left(  \xi\right)  \widetilde{x}_{m+1}\left(
\omega\right)  \right\}
\]
which extends the orthogonality relation (\ref{orthogonality_eigenvectors_P}).
Transformed back to $x_{j}\left(  \xi\right)  $ using (\ref{def_h_i}) yields%
\begin{equation}
\left(  \xi-\omega\right)  \sum_{j=0}^{m}h_{j+1}^{2}x_{j}\left(  \xi\right)
x_{j}\left(  \omega\right)  =p_{m}h_{m+1}^{2}\left\{  x_{m}\left(
\omega\right)  x_{m+1}\left(  \xi\right)  -x_{m}\left(  \xi\right)
x_{m+1}\left(  \omega\right)  \right\}
\label{Christoffel-Darboux_formula_eigenvectors_P}%
\end{equation}
Since $\omega=0$ is an eigenvalue with corresponding eigenvector
$x(0)=\frac{1}{N+1}u$, each other real eigenvalue $\xi\neq0$ must obey%
\[
\xi\sum_{j=0}^{m}h_{j+1}^{2}x_{j}\left(  \xi\right)  =p_{m}h_{m+1}^{2}\left\{
x_{m+1}\left(  \xi\right)  -x_{m}\left(  \xi\right)  \right\}
\]
Taking $p_{N}=0$ into account, the Christoffel-Darboux formula
(\ref{Christoffel-Darboux_formula_eigenvectors_P}) extends
(\ref{orthogonality_eigenvectors_P}) to all $0\leq m\leq N$.

\section{Second largest zero of $c_{P}\left(  \xi\right)  $}

\label{sec_properties_zeros_characteristic_coefficient_P}

The zero of a complex function can be expressed as a Lagrange series
\cite{Whittaker_Watson,Markushevich}. When all Taylor coefficients $f_{k}$ of
a function expanded around a point $z_{0}$ are known, our framework of
characteristic coefficients, first published in \cite{PVM_ASYM}, provides all
coefficients in the corresponding Lagrange series in terms of $f_{k}$. In
particular, the second largest zero $\zeta$ closest to $\xi=0$, based on the
Lagrange expansion (see e.g. \cite[p. 305]{PVM_graphspectra}) up to order 4 in
$\frac{f_{0}}{f_{1}}$, is
\begin{equation}
\zeta\approx-\frac{f_{0}}{f_{1}}-\frac{f_{2}}{f_{1}}\;\left(  \frac{f_{0}%
}{f_{1}}\right)  ^{2}+\left[  -2\,\left(  \frac{f_{2}}{f_{1}}\right)
^{2}+\frac{f_{3}}{f_{1}}\right]  \;\left(  \frac{f_{0}}{f_{1}}\right)
^{3}+O\left(  \left(  \frac{f_{0}}{f_{1}}\right)  ^{4}\right)
\label{zero_Lagrange_series_up_to_4_order}%
\end{equation}
Since all Taylor coefficients $f_{k}$ of the characteristic polynomial
$c_{P}\left(  \xi\right)  $ around $\xi=0$ are known, we can formally compute
the zero $\zeta$ to \emph{any order or accuracy}. The fact that all
coefficients $f_{k}$ are non-zero and that $\xi=0$ is the largest zero of
$c_{P}\left(  \xi\right)  $ guarantees that the Lagrange series converges
fast. In fact, the first term in (\ref{zero_Lagrange_series_up_to_4_order})
equals the first iteration in the Newton-Raphson method and the point
$z_{0}=0$ is an ideal expansion point. This article demonstrates this
computation up to second order, hence, using the explicit knowledge of
$f_{0},f_{1}$ and $f_{2}$. Proceeding further with $f_{3}$ is possible,
however, at the expense of huge computations, from which we refrained, mainly
because numerical computations in Section \ref{sec_epsilon_SIS_epidemics}
demonstrate a good accuracy of $\zeta$ only based on the three coefficients
$f_{0},f_{1}$ and $f_{2}$.

The sum\footnote{The sum of the zeros of $\frac{c_{P}\left(  \xi\right)  }%
{\xi}$ (taking into account that $p_{N}=0$) equals%
\[
\sum_{k=1}^{N}z_{k}=-\frac{f_{N-1}}{f_{N}}=-\sum_{m=0}^{N}\left(  p_{m}%
+q_{m}\right)
\]
Since $0\leq p_{m}+q_{m}=1-r_{m}\leq1$ and $p_{N}=q_{0}=0$, the average of the
zeros lies between zero and minus one. The product of the zeros follows from
(\ref{characteristic_polynomial_P_grw_in_zeros}) as%
\[
\prod\limits_{k=1}^{N}\left(  -z_{k}\right)  =\frac{f_{0}}{f_{N}}=\frac
{\prod_{m=0}^{N-1}q_{m+1}}{\pi_{0}}=\prod_{m=0}^{N-1}q_{m+1}+\sum_{j=1}%
^{N-1}\prod_{m=j}^{N-1}q_{m+1}\prod_{m=0}^{j-1}p_{m}+\prod_{m=0}^{N-1}p_{m}%
\]
which is, by the Perron-Frobenius Theorem strictly smaller than 1. Finally, we
also compute $\sum_{k=1}^{N}z_{k}^{2}=\left(  \frac{f_{N-1}}{f_{N}}\right)
^{2}-2\frac{f_{N-2}}{f_{N}}$ from the Newton identities with
(\ref{c_k(j)_for_k=j-2}) as%
\[
\sum_{k=1}^{N}z_{k}^{2}=\sum_{m=0}^{N}\left\{  \left(  q_{m}+p_{m}\right)
^{2}+2q_{m}p_{m-1}\right\}
\]
} of the inverse of the zeros of $\frac{c_{P}\left(  \xi\right)  }{\xi}$
follows from the Newton identities \cite[p. 305]{PVM_graphspectra} as%
\[
\sum_{k=1}^{N}\frac{1}{z_{k}}=-\frac{f_{1}}{f_{0}}%
\]
from which%
\[
-\zeta=\frac{1}{\frac{f_{1}}{f_{0}}+\sum_{k=2}^{N}\frac{1}{z_{k}}}%
\]
Since all zeros $z_{k}$ of $\frac{c_{P}\left(  \xi\right)  }{\xi}$ are
negative, we have%
\[
-\zeta=\frac{1}{\frac{f_{1}}{f_{0}}+\sum_{k=2}^{N}\frac{1}{z_{k}}}>\frac
{f_{0}}{f_{1}}%
\]
so that
\begin{equation}
\zeta<-\frac{f_{0}}{f_{1}} \label{upper_bound_zeta_first_order}%
\end{equation}
demonstrating that $-\frac{f_{0}}{f_{1}}$ is an upper bound for $\zeta$. This
observation also follows from the above Lagrange series
(\ref{zero_Lagrange_series_up_to_4_order}).

From (\ref{def_f_k}), we have that $f_{0}=\frac{1}{\pi_{0}}$, where $\pi_{0}$
is the zero component of the state-state vector of $P$ (eigenvector belonging
to eigenvalue $\lambda=1$), and
\[
f_{1}=\frac{1}{q_{1}}+\sum_{j=2}^{N}\frac{c_{1}(j)}{\prod_{m=0}^{j-1}q_{m+1}}%
\]
which becomes with (\ref{c1(j)_explicit}),%
\begin{equation}
f_{1}=\frac{1}{q_{1}}+\sum_{j=2}^{N}\frac{1}{\prod_{m=0}^{j-1}q_{m+1}}%
\sum_{r=0}^{j-1}\prod\limits_{s=j-r}^{j-1}p_{s}\sum_{k=0}^{j-1-r}%
\prod\limits_{m=0}^{j-2-r-k}p_{m}\prod\limits_{l=j-r-k}^{j-1-r}q_{l}
\label{f1}%
\end{equation}
The number of terms in $f_{1}$ equals $1+\sum_{j=1}^{N-1}\sum_{r=0}^{j}%
\sum_{k=0}^{j-r}1=\frac{N\left(  N+1\right)  \left(  N+2\right)  }{6}%
=\binom{N+2}{3}$. Hence, the \emph{lower} bound for $-\zeta$ is%
\[
\frac{f_{0}}{f_{1}}=\frac{1+\sum_{j=1}^{N}\prod_{m=0}^{j-1}\frac{p_{m}%
}{q_{m+1}}}{\frac{1}{q_{1}}+\sum_{j=2}^{N}\frac{1}{\prod_{m=0}^{j-1}q_{m+1}%
}\sum_{r=0}^{j-1}\prod\limits_{s=j-r}^{j-1}p_{s}\sum_{k=0}^{j-1-r}%
\prod\limits_{m=0}^{j-2-r-k}p_{m}\prod\limits_{l=j-r-k}^{j-1-r}q_{l}}%
\]

Similarly, combining (\ref{def_f_k}), (\ref{c2(j)_explicit}) and (\ref{f1})
yields $\frac{f_{2}}{f_{1}}$, and so establishing the Lagrange series for
$\zeta$ up to second order in $\frac{f_{0}}{f_{1}}$,
\begin{equation}
\zeta\approx-\frac{f_{0}}{f_{1}}-\frac{f_{2}}{f_{1}}\;\left(  \frac{f_{0}%
}{f_{1}}\right)  ^{2} \label{zero_Lagrange_series_up_to_2_order}%
\end{equation}
The inverse of the squares of the zeros of $\frac{c_{P}\left(  \xi\right)
}{\xi}$ equals \cite[p. 305]{PVM_graphspectra}
\[
\sum_{k=1}^{N}\frac{1}{z_{k}^{2}}=\left(  \frac{f_{1}}{f_{0}}\right)
^{2}-\frac{2f_{2}}{f_{0}}\geq0
\]
from which%
\[
\frac{1}{\zeta^{2}}=\left(  \frac{f_{1}}{f_{0}}\right)  ^{2}-\frac{2f_{2}%
}{f_{0}}-\sum_{k=2}^{N}\frac{1}{z_{k}^{2}}%
\]
and%
\[
\zeta^{2}=\frac{1}{\left(  \frac{f_{1}}{f_{0}}\right)  ^{2}-\frac{2f_{2}%
}{f_{0}}-\sum_{k=2}^{N}\frac{1}{z_{k}^{2}}}\geq\left(  \frac{f_{0}}{f_{1}%
}\right)  ^{2}\frac{1}{1-\frac{2f_{2}}{f_{0}}\left(  \frac{f_{0}}{f_{1}%
}\right)  ^{2}}%
\]
where the inequality follows because all zeros are real. Thus, a sharper upper
bound for $\zeta$ is found%
\begin{equation}
\zeta\leq-\frac{f_{0}}{f_{1}}\frac{1}{\sqrt{1-\frac{2f_{2}}{f_{0}}\left(
\frac{f_{0}}{f_{1}}\right)  ^{2}}} \label{upper_bound_Newton_zeta}%
\end{equation}
After expansion of the right-hand side in (\ref{upper_bound_Newton_zeta}), we
find%
\[
\left\vert \zeta\right\vert \geq\left(  \frac{f_{0}}{f_{1}}\right)  \left(
1+\frac{f_{2}}{f_{0}}\left(  \frac{f_{0}}{f_{1}}\right)  ^{2}+\frac{3}%
{2}\left\{  \frac{f_{2}}{f_{0}}\left(  \frac{f_{0}}{f_{1}}\right)
^{2}\right\}  ^{2}+O\left(  \left\{  \frac{f_{2}}{f_{0}}\left(  \frac{f_{0}%
}{f_{1}}\right)  ^{2}\right\}  ^{3}\right)  \right)
\]
which should be compared with the Lagrange expansion up to third order,%
\begin{align*}
-\zeta &  \approx\left(  \frac{f_{0}}{f_{1}}\right)  \left(  1+\frac{f_{2}%
}{f_{1}}\;\left(  \frac{f_{0}}{f_{1}}\right)  +\left[  2\,\left(  \frac{f_{2}%
}{f_{1}}\right)  ^{2}-\frac{f_{3}}{f_{1}}\right]  \;\left(  \frac{f_{0}}%
{f_{1}}\right)  ^{2}\right) \\
&  =\left(  \frac{f_{0}}{f_{1}}\right)  \left(  1+\frac{f_{2}}{f_{0}}\;\left(
\frac{f_{0}}{f_{1}}\right)  ^{2}+2\,\left\{  \frac{f_{2}}{f_{0}}\left(
\frac{f_{0}}{f_{1}}\right)  ^{2}\right\}  ^{2}-\frac{f_{3}}{f_{1}}\left(
\frac{f_{0}}{f_{1}}\right)  ^{2}\right)
\end{align*}
In summary, based on the knowledge of the coefficients $f_{0}$, $f_{1}$ and
$f_{2}$, the second largest zero $\zeta$ of $c_{P}\left(  \xi\right)  $ is
approximated by a Lagrange series (\ref{zero_Lagrange_series_up_to_2_order})
up to second order, possesses an upper\footnote{The interlacing property of
the orthogonal polynomials provides us with lower bounds for $\zeta$. The
interlacing theorem for orthogonal polynomials states that between two zeros
of $\rho_{k}\left(  \xi\right)  $, defined in
(\ref{def_rj_orthogonal_polynomial_of_P}), there is at least on zero of
$\rho_{l}\left(  \xi\right)  $ with $l>k$. The best lower bound for $\zeta$
thus equals the second largest zero of $\rho_{N}\left(  \xi\right)  $, which
is, unfortunately, more difficult to compute than $\zeta$ itself, because the
largest zero is negative and unknown in contrast to $\rho_{N+1}\left(
\xi\right)  $ where it is zero.} bound (\ref{upper_bound_zeta_first_order})
and a sharper bound (\ref{upper_bound_Newton_zeta}).

Incidentally, we have also shown how subsequent terms in the Lagrange series
can be computed from the Newton identities for the sum of inverse powers of
the zeros. The combination of the knowledge of the Newton identities with the
Lagrange series around a certain complex number can shed additional insight
into the convergence of the Lagrange series.

\section{The function $F\left(  \tau\right)  $}

\label{sec_F(tau)}We have shown that $\zeta\approx-\frac{1}{f_{1}}$ for
$\tau>\frac{1}{N}$, which led to the result
(\ref{zeta_SIS_KN_above_threshold_x>1}). We reconsider $F\left(  \tau\right)
=\lim_{\varepsilon\rightarrow0}f_{1}$ in (\ref{f1_eps=0}), which is also
rewritten as%
\begin{equation}
\beta F\left(  \tau\right)  =\sum_{j=1}^{N}\sum_{r=1}^{j}\frac{\left(
N-1-j+r\right)  !}{j\left(  N-j\right)  !}\tau^{r} \label{def_beta_F}%
\end{equation}
In this section, we explore properties of $F\left(  \tau\right)  $: in Section
\ref{sec_Laplace_transform_F(tau)}, $F\left(  \tau\right)  $ is expressed as a
Laplace transform, from which alternative exact forms for $F\left(
\tau\right)  $ are deduced in Section
\ref{sec_F(tau)_as_exponential_integrals}. Finally, Section
\ref{sec_order_F(tau_largeN} presents the asymptotic form
(\ref{asymptotic_F(tau)_Nlarge_x_fixed}) of $F\left(  \tau\right)  $ for
$\tau=\frac{x}{N}$, with fixed $x$ and large $N$.

\subsection{$F\left(  \tau\right)  $ as a Laplace transform}

\label{sec_Laplace_transform_F(tau)}

\begin{theorem}
\label{theorem_Laplace_f1_eps=0} For $\tau>\tau_{c}$, $F\left(  \tau\right)
=\lim_{\varepsilon\rightarrow0}f_{1}$ can be expressed as a Laplace transform
\begin{equation}
F\left(  \tau\right)  =\frac{1}{\beta}\int_{0}^{\infty}\left(  \int
_{1}^{\infty}\frac{dx}{x^{N+1}}\frac{\left(  u+x\right)  ^{N}-1}{\left(
u+x\right)  -1}\right)  e^{-\frac{1}{\tau}u}du \label{Laplace_transform_f1}%
\end{equation}

\end{theorem}

\textbf{Proof}: We rewrite (\ref{f1_eps=0}) as%
\[
F\left(  \tau\right)  =\frac{1}{\delta}\sum_{j=1}^{N}\frac{x_{j}}{j}%
\]
with%
\[
x_{j}=\sum_{r=0}^{j-1}\frac{\left(  N-j+r\right)  !}{\left(  N-j\right)
!}\tau^{r}%
\]
which obeys the recursion%
\begin{equation}
x_{j+1}=x_{j}(N-j)\tau+1 \label{recursion_x_j}%
\end{equation}
with $x_{1}=1$. Furthermore, $x_{k}=0$ for $k<1$. Indeed,%
\begin{align*}
x_{j+1}  &  =\sum_{r=0}^{j}\frac{\left(  N-j-1+r\right)  !}{\left(
N-j-1\right)  !}\tau^{r}=\tau\left(  N-j\right)  \sum_{r=-1}^{j-1}%
\frac{\left(  N-j+r\right)  !}{\left(  N-j\right)  !}\tau^{r}\\
&  =x_{j}(N-j)\tau+\tau\left(  N-j\right)  \frac{\left(  N-j-1\right)
!}{\left(  N-j\right)  !}\tau^{-1}=x_{j}(N-j)\tau+1
\end{align*}

Now, consider%
\begin{equation}
g\left(  z,N\right)  =\sum_{j=0}^{N-1}x_{j+1}z^{j}=\frac{1}{z}\sum_{j=1}%
^{N}x_{j}z^{j} \label{def_g_z_N}%
\end{equation}
with $g\left(  0,N\right)  =x_{1}=1$ so that%
\[
\int_{0}^{1}g\left(  z,N\right)  dz=\sum_{j=0}^{N-1}\frac{x_{j+1}}{j+1}%
=\sum_{j=1}^{N}\frac{x_{j}}{j}=\delta F\left(  \tau\right)
\]
After multiplying both sides of the recursion (\ref{recursion_x_j}) with
$z^{j}$ and summing over $j$, we obtain%
\begin{align*}
g\left(  z,N\right)   &  =\sum_{j=0}^{N-1}x_{j+1}z^{j}=\tau\sum_{j=0}%
^{N-1}x_{j}(N-j)z^{j}+\sum_{j=0}^{N-1}z^{j}\\
&  =\tau N\sum_{j=1}^{N}x_{j}z^{j}-\tau\sum_{j=0}^{N}jx_{j}z^{j}+\frac
{z^{N}-1}{z-1}%
\end{align*}
since $x_{0}=0$. With the definition (\ref{def_g_z_N}) and%
\[
\frac{d}{dz}\left(  zg\left(  z,N\right)  \right)  =\frac{1}{z}\sum_{j=0}%
^{N}jx_{j}z^{j}%
\]
we find%
\[
g\left(  z,N\right)  =\tau Nzg\left(  z,N\right)  -\tau z\frac{d}{dz}\left(
zg\left(  z,N\right)  \right)  +\frac{z^{N}-1}{z-1}%
\]
Thus,%
\[
\frac{d}{dz}\left(  zg\left(  z,N\right)  \right)  =\left(  N-\frac{1}{\tau
z}\right)  g\left(  z,N\right)  +\frac{z^{N}-1}{\tau z\left(  z-1\right)  }%
\]
The differential equation for $g$ becomes%
\[
-\tau z^{2}\frac{d}{dz}\left(  g\left(  z,N\right)  \right)  +\left(
\tau\left(  N-1\right)  z-1\right)  g\left(  z,N\right)  =-\frac{z^{N}-1}{z-1}%
\]
The homogeneous differential equation is rewritten as%
\[
\frac{d}{dz}\left(  \log h\left(  z,N\right)  \right)  =\frac{\tau\left(
N-1\right)  z-1}{\tau z^{2}}=\frac{\left(  N-1\right)  }{z}-\frac{1}{\tau
z^{2}}%
\]
Integration yields%
\[
\log h\left(  z,N\right)  =\left(  N-1\right)  \ln z+\frac{1}{\tau z}+\ln C
\]
where $C$ is a constant. Thus,%
\[
h\left(  z,N\right)  =Cz^{N-1}e^{\frac{1}{\tau z}}%
\]
Using the variation of a constant method yields%
\[
g\left(  z,N\right)  =C\left(  z\right)  z^{N-1}e^{\frac{1}{\tau z}}%
\]
where the function $C\left(  z\right)  $ must obey the differential equation.
Hence, after substitution, we have%
\[
g^{\prime}\left(  z,N\right)  =C^{\prime}\left(  z\right)  z^{N-1}e^{\frac
{1}{\tau z}}+C\left(  z\right)  \frac{d}{dz}\left(  z^{N-1}e^{\frac{1}{\tau
z}}\right)
\]
and%
\[
-\tau z^{2}C^{\prime}\left(  z\right)  z^{N-1}e^{\frac{1}{\tau z}}-\tau
z^{2}C\left(  z\right)  \frac{d}{dz}\left(  z^{N-1}e^{\frac{1}{\tau z}%
}\right)  +\left(  \tau\left(  N-1\right)  z-1\right)  C\left(  z\right)
z^{N-1}e^{\frac{1}{\tau z}}=-\frac{z^{N}-1}{z-1}%
\]
Since
\begin{align*}
-\tau z^{2}C\left(  z\right)  \frac{d}{dz}\left(  z^{N-1}e^{\frac{1}{\tau z}%
}\right)   &  =-\tau z^{2}C\left(  z\right)  \left(  N-1\right)
z^{N-2}e^{\frac{1}{\tau z}}+\tau z^{2}C\left(  z\right)  z^{N-1}e^{\frac
{1}{\tau z}}\left(  \frac{1}{\tau z^{2}}\right) \\
&  =\left\{  -\tau z\left(  N-1\right)  +1\right\}  z^{N-1}C\left(  z\right)
e^{\frac{1}{\tau z}}%
\end{align*}
we see (as required by the method of the variation of a constant) that%
\[
C^{\prime}\left(  z\right)  =\frac{1}{\tau}\frac{z^{N}-1}{\left(  z-1\right)
z^{N+1}}e^{-\frac{1}{\tau z}}%
\]
from which $C^{\prime}\left(  0\right)  =0$ for any finite $N$. After
integration, we arrive at%
\[
g\left(  z,N\right)  =\frac{1}{\tau}z^{N-1}e^{\frac{1}{\tau z}}\left(
\int_{0}^{z}\frac{u^{N}-1}{\left(  u-1\right)  u^{N+1}}e^{-\frac{1}{\tau u}%
}du+C\left(  0\right)  \right)
\]
where the constant $C\left(  0\right)  $ needs to be chosen so that $g\left(
0,N\right)  =1$. The only possible value is $C\left(  0\right)  =0$ in order
to have a finite limit for $\lim_{z\rightarrow0}g\left(  z,N\right)  $. Then,%
\[
g\left(  z,N\right)  =\frac{1}{\tau}z^{N-1}e^{\frac{1}{\tau z}}\int_{0}%
^{z}\frac{u^{N}-1}{\left(  u-1\right)  u^{N+1}}e^{-\frac{1}{\tau u}}du
\]
Substituting $u=\frac{1}{y}$,
\begin{equation}
g\left(  z,N\right)  =\frac{1}{\tau}z^{N-1}e^{\frac{1}{\tau z}}\int_{\frac
{1}{z}}^{\infty}\frac{\left(  1-y^{N}\right)  }{\left(  1-y\right)  }%
e^{-\frac{1}{\tau}y}dy \label{g_integraal}%
\end{equation}
Finally,%
\begin{align*}
\delta F\left(  \tau\right)   &  =\frac{1}{\tau}\int_{0}^{1}dz\ z^{N-1}%
e^{\frac{1}{\tau z}}\int_{\frac{1}{z}}^{\infty}\frac{\left(  1-y^{N}\right)
}{\left(  1-y\right)  }e^{-\frac{1}{\tau}y}dy\\
&  =\frac{1}{\tau}\int_{1}^{\infty}dx\ \frac{1}{x^{N+1}}\int_{x}^{\infty}%
\frac{\left(  1-y^{N}\right)  }{\left(  1-y\right)  }e^{-\frac{1}{\tau}\left(
y-x\right)  }dy
\end{align*}
Let $u=y-x$ in the $y$-integral, then%
\[
\delta F\left(  \tau\right)  =\frac{1}{\tau}\int_{1}^{\infty}dx\ \frac
{1}{x^{N+1}}\int_{0}^{\infty}\frac{\left(  1-\left(  u+x\right)  ^{N}\right)
}{\left(  1-u-x\right)  }e^{-\frac{1}{\tau}u}du
\]
so that $F\left(  \tau\right)  $ can be written as a Laplace transform
(\ref{Laplace_transform_f1}).\hfill$\square\medskip$

\subsection{$F\left(  \tau\right)  $ in terms of exponential integrals}

\label{sec_F(tau)_as_exponential_integrals}The integrand in
(\ref{Laplace_transform_f1}) can be rewritten as%
\[
h_{N}\left(  u\right)  =\int_{1}^{\infty}\frac{dx}{x^{N+1}}\frac{\left(
u+x\right)  ^{N}-1}{\left(  u+x\right)  -1}=\int_{0}^{1}\frac{\left(
yu+1\right)  ^{N}-y^{N}}{\left(  yu+1\right)  -y}dy
\]
so that we obtain an alternative integral%
\[
\beta F\left(  \tau\right)  =\int_{0}^{\infty}\left(  \int_{0}^{1}%
\frac{\left(  yu+1\right)  ^{N}-y^{N}}{\left(  yu+1\right)  -y}dy\right)
e^{-\frac{1}{\tau}u}du
\]
or%
\begin{equation}
\beta F\left(  \tau\right)  =\int_{0}^{1}\left(  \int_{0}^{\infty}%
\frac{\left(  yu+1\right)  ^{N}e^{-\frac{1}{\tau}u}}{yu+1-y}du\right)
dy-\int_{0}^{1}y^{N}\left(  \int_{0}^{\infty}\frac{e^{-\frac{1}{\tau}u}%
}{yu+1-y}du\right)  dy \label{integral_betaF}%
\end{equation}
which will be exploited below.

In the next Theorem \ref{theorem_betaF_series_exponential_integrals}, we show
that $F\left(  \tau\right)  $ can be expressed in terms of the exponential
integrals $E_{n}\left(  x\right)  $ of integer order $n$, defined in
\cite[Chapter 5]{Abramowitz} as%
\[
E_{n}\left(  x\right)  =\int_{1}^{\infty}\frac{e^{-xt}}{t^{n}}dt
\]

\begin{theorem}
\label{theorem_betaF_series_exponential_integrals} For $\tau>\frac{1}{N}$,
$F\left(  \tau\right)  =\lim_{\varepsilon\rightarrow0}f_{1}$ equals%
\begin{equation}
F\left(  \tau\right)  =\frac{N!}{\beta}\sum_{k=1}^{N+1}\frac{L_{k}\left(
\tau\right)  }{\left(  N+1-k\right)  !}-\frac{1}{\beta}\int_{0}^{\infty}%
\frac{e^{w}E_{N+1}\left(  w\right)  dw}{w+\frac{1}{\tau}}
\label{betaF_in_exponential_integrals}%
\end{equation}
where%
\begin{equation}
L_{k}\left(  \tau\right)  =\int_{0}^{\infty}\frac{e^{w}E_{k}\left(  w\right)
}{\left(  w+\frac{1}{^{\tau}}\right)  ^{k}}dw \label{def_Lk(tau)}%
\end{equation}

\end{theorem}

\textbf{Proof}: We start with the second $u$-integral in (\ref{integral_betaF}%
) for $\beta F\left(  \tau\right)  $,%
\begin{align*}
\int_{0}^{\infty}\frac{e^{-\frac{1}{\tau}u}}{yu+1-y}du  &  =\frac{e^{\frac
{1}{\tau}\left(  \frac{1}{y}-1\right)  }}{y}\int_{0}^{\infty}\frac
{e^{-\frac{1}{\tau}\left(  u+\frac{1}{y}-1\right)  }}{u+\frac{1}{y}-1}du\\
&  =\frac{e^{\frac{1}{\tau}\left(  \frac{1}{y}-1\right)  }}{y}\int_{0}%
^{\infty}\int_{\frac{1}{\tau}}^{\infty}e^{-p\left(  u+\frac{1}{y}-1\right)
}dpdu\\
&  =\frac{e^{\frac{1}{\tau}\left(  \frac{1}{y}-1\right)  }}{y}\int_{\frac
{1}{\tau}}^{\infty}e^{-p\left(  \frac{1}{y}-1\right)  }\int_{0}^{\infty
}e^{-pu}dudp\\
&  =\frac{1}{y}\int_{\frac{1}{\tau}}^{\infty}\frac{e^{-\left(  p-\frac{1}%
{\tau}\right)  \left(  \frac{1}{y}-1\right)  }}{p}dp\\
&  =\frac{1}{y}\int_{0}^{\infty}\frac{e^{-w\left(  \frac{1}{y}-1\right)  }%
}{w+\frac{1}{\tau}}dw
\end{align*}
where the reversal of integrations is allowed by absolute convergence
(Fubinni's Theorem). The second $u$-integral in (\ref{integral_betaF}) for
$\beta F\left(  \tau\right)  $ becomes
\begin{align*}
\int_{0}^{1}y^{N}\left(  \int_{0}^{\infty}\frac{e^{-\frac{1}{\tau}u}}%
{yu+1-y}du\right)  dy  &  =\int_{0}^{\infty}\frac{e^{w}dw}{w+\frac{1}{\tau}%
}\left(  \int_{0}^{1}e^{-w\frac{1}{y}}y^{N-1}dy\right) \\
&  =\int_{0}^{\infty}\frac{e^{w}dw}{w+\frac{1}{\tau}}\left(  \int_{1}^{\infty
}\frac{e^{-wx}}{x^{N+1}}dx\right)
\end{align*}
or%
\begin{equation}
\int_{0}^{1}y^{N}\left(  \int_{0}^{\infty}\frac{e^{-\frac{1}{\tau}u}}%
{yu+1-y}du\right)  dy=\int_{0}^{\infty}\frac{e^{w}E_{N+1}\left(  w\right)
dw}{w+\frac{1}{\tau}} \label{tweede_stuk_beta_F(tau)_intergral}%
\end{equation}

We now focus on the first integral in (\ref{integral_betaF}) and start with%
\begin{align*}
\int_{0}^{\infty}\frac{\left(  yu+1\right)  ^{N}e^{-\frac{1}{\tau}u}}%
{yu+1-y}du  &  =\frac{e^{\frac{1}{\tau}\left(  \frac{1}{y}-1\right)  }}{y}%
\int_{0}^{\infty}\frac{\left(  yu+1\right)  ^{N}e^{-\frac{1}{\tau}\left(
u+\frac{1}{y}-1\right)  }}{u+\frac{1}{y}-1}du\\
&  =\frac{e^{\frac{1}{\tau}\left(  \frac{1}{y}-1\right)  }}{y}\int_{0}%
^{\infty}\left(  yu+1\right)  ^{N}\int_{\frac{1}{\tau}}^{\infty}e^{-p\left(
u+\frac{1}{y}-1\right)  }dpdu\\
&  =\frac{e^{\frac{1}{\tau}\left(  \frac{1}{y}-1\right)  }y^{N}}{y}\int
_{\frac{1}{\tau}}^{\infty}e^{-p\left(  \frac{1}{y}-1\right)  }\left(  \int
_{0}^{\infty}\left(  u+\frac{1}{y}\right)  ^{N}e^{-pu}du\right)  dp
\end{align*}
and%
\begin{align*}
\int_{0}^{\infty}\left(  u+\frac{1}{y}\right)  ^{N}e^{-pu}du  &  =e^{\frac
{p}{y}}\int_{\frac{1}{y}}^{\infty}w^{N}e^{-pw}dw=\frac{e^{\frac{p}{y}}%
}{p^{N+1}}\int_{\frac{p}{y}}^{\infty}t^{N}e^{-t}dt\\
&  =\frac{N!}{p^{N+1}}\sum_{j=0}^{N}\frac{\left(  \frac{p}{y}\right)  ^{j}%
}{j!}%
\end{align*}
where the incomplete Gamma function for integer $k$%
\[
\int_{x}^{\infty}t^{k}e^{-t}dt=k!e^{-x}\sum_{j=0}^{k}\frac{x^{j}}{j!}%
\]
is used. Hence\footnote{Notice that we cannot use the approximation for large
$N$,%
\[
\sum_{j=0}^{N}\frac{\left(  \frac{p}{y}\right)  ^{j}}{j!}\lesssim e^{\frac
{p}{y}}%
\]
because then the $p$-integral diverges.},%
\begin{align*}
\int_{0}^{\infty}\frac{\left(  yu+1\right)  ^{N}e^{-\frac{1}{\tau}u}}%
{yu+1-y}du  &  =\frac{e^{\frac{1}{\tau}\left(  \frac{1}{y}-1\right)  }y^{N}%
}{y}\int_{\frac{1}{\tau}}^{\infty}e^{-p\left(  \frac{1}{y}-1\right)  }\left(
\frac{N!}{p^{N+1}}\sum_{j=0}^{N}\frac{\left(  \frac{p}{y}\right)  ^{j}}%
{j!}\right)  dp\\
&  =N!e^{\frac{1}{\tau}\left(  \frac{1}{y}-1\right)  }y^{N-1}\sum_{j=0}%
^{N}\frac{\left(  \frac{1}{y}\right)  ^{j}}{j!}\int_{\frac{1}{\tau}}^{\infty
}\frac{e^{-p\left(  \frac{1}{y}-1\right)  }}{p^{N+1-j}}dp
\end{align*}

We consider now the first integral in the expression (\ref{integral_betaF})
for $\beta F\left(  \tau\right)  $,%
\begin{align*}
\int_{0}^{1}\left(  \int_{0}^{\infty}\frac{\left(  yu+1\right)  ^{N}%
e^{-\frac{1}{\tau}u}}{yu+1-y}du\right)  dy  &  =N!\int_{0}^{1}e^{\frac{1}%
{\tau}\left(  \frac{1}{y}-1\right)  }y^{N-1}\sum_{j=0}^{N}\frac{\left(
\frac{1}{y}\right)  ^{j}}{j!}\int_{\frac{1}{\tau}}^{\infty}\frac{e^{-p\left(
\frac{1}{y}-1\right)  }}{p^{N+1-j}}dpdy\\
&  =N!\sum_{j=0}^{N}\frac{1}{j!}\int_{\frac{1}{\tau}}^{\infty}\frac
{e^{p-\frac{1}{\tau}}}{p^{N+1-j}}dp\int_{0}^{1}e^{\left(  \frac{1}{\tau
}-p\right)  \frac{1}{y}}y^{N-1-j}dy\\
&  =N!\sum_{j=0}^{N}\frac{1}{j!}\int_{0}^{\infty}\frac{e^{w}}{\left(
w+\frac{1}{^{\tau}}\right)  ^{N+1-j}}dw\int_{0}^{1}e^{-w\frac{1}{y}}%
y^{N-1-j}dy
\end{align*}
Since%
\[
\int_{0}^{1}e^{-w\frac{1}{y}}y^{N-1-j}dy=\int_{1}^{\infty}\frac{e^{-wu}%
}{u^{N+1-j}}du=E_{N+1-j}\left(  w\right)
\]
we find that%
\begin{align}
\int_{0}^{1}\left(  \int_{0}^{\infty}\frac{\left(  yu+1\right)  ^{N}%
e^{-\frac{1}{\tau}u}}{yu+1-y}du\right)  dy  &  =N!\sum_{j=0}^{N}\frac{1}%
{j!}\int_{0}^{\infty}\frac{e^{w}E_{N+1-j}\left(  w\right)  }{\left(
w+\frac{1}{^{\tau}}\right)  ^{N+1-j}}dw\nonumber\\
&  =N!\sum_{k=1}^{N+1}\frac{L_{k}\left(  \tau\right)  }{\left(  N+1-k\right)
!} \label{dominant_betaF_integral_expanded_Lk(tau)}%
\end{align}
which proves the theorem.\hfill$\square\medskip$

An expression that avoids the summation is%
\[
\int_{0}^{\infty}\frac{\left(  yu+1\right)  ^{N}e^{-\frac{1}{\tau}u}}%
{yu+1-y}du=\frac{e^{\frac{1}{\tau}\left(  \frac{1}{y}-1\right)  }y^{N}}{y}%
\int_{\frac{1}{\tau}}^{\infty}dp\frac{e^{p}}{p^{N+1}}\int_{\frac{p}{y}%
}^{\infty}t^{N}e^{-t}dt
\]
Then,%
\begin{align*}
\int_{0}^{1}\left(  \int_{0}^{\infty}\frac{\left(  yu+1\right)  ^{N}%
e^{-\frac{1}{\tau}u}}{yu+1-y}du\right)  dy  &  =e^{-\frac{1}{\tau}}\int
_{0}^{1}dye^{\frac{1}{\tau y}}y^{N-1}\int_{\frac{1}{\tau}}^{\infty}%
dp\frac{e^{p}}{p^{N+1}}\int_{\frac{p}{y}}^{\infty}t^{N}e^{-t}dt\\
&  =e^{-\frac{1}{\tau}}\int_{\frac{1}{\tau}}^{\infty}dp\frac{e^{p}}{p^{N+1}%
}\int_{0}^{1}dye^{\frac{1}{\tau y}}y^{N-1}\int_{\frac{p}{y}}^{\infty}%
t^{N}e^{-t}dt\\
&  =\frac{e^{-\frac{1}{\tau}}}{\tau^{N}}\int_{\frac{1}{\tau}}^{\infty}%
dp\frac{e^{p}}{p^{N+1}}\int_{\frac{1}{\tau}}^{\infty}du\frac{e^{u}}{u^{N+1}%
}\int_{\tau pu}^{\infty}t^{N}e^{-t}dt
\end{align*}
We will not further use this triple integral, although it suggests a change of
variables $s=p+u$ and $r=pu$, which, as we found, did not lead to useful results.

\subsubsection{Other exact series for $\beta F\left(  \tau\right)  $}

The expression (\ref{betaF_in_exponential_integrals}) for $\beta F\left(
\tau\right)  $ in Theorem \ref{theorem_betaF_series_exponential_integrals}
will be further explored by using properties of the exponential integral.

After partial integration of $E_{k}\left(  w\right)  =\int_{1}^{\infty}%
\frac{e^{-wt}}{t^{k}}dt$, we obtain the recursion,%
\begin{equation}
E_{k+1}\left(  w\right)  =\frac{1}{k}e^{-w}-\frac{w}{k}E_{k}\left(  w\right)
\label{recursion_exponential_integral}%
\end{equation}
After $p$ iteration, we find%
\begin{equation}
E_{k}\left(  w\right)  =e^{-w}\sum_{j=0}^{p-1}\frac{\left(  -1\right)
^{j}\left(  k-j-2\right)  !w^{j}}{\left(  k-1\right)  !}+\frac{\left(
-1\right)  ^{p}\left(  k-p-1\right)  !w^{p}}{\left(  k-1\right)  !}%
E_{k-p}\left(  w\right)  \label{iterated_p_recursion_exponential_integral}%
\end{equation}
Introducing (\ref{iterated_p_recursion_exponential_integral}) into
(\ref{def_Lk(tau)}) gives%
\[
L_{k}\left(  \tau\right)  =\sum_{j=0}^{p-1}\frac{\left(  -1\right)
^{j}\left(  k-j-2\right)  !}{\left(  k-1\right)  !}\int_{0}^{\infty}%
\frac{w^{j}dw}{\left(  w+\frac{1}{^{\tau}}\right)  ^{k}}+\frac{\left(
-1\right)  ^{p}\left(  k-p-1\right)  !}{\left(  k-1\right)  !}\int_{0}%
^{\infty}\frac{w^{p}e^{w}}{\left(  w+\frac{1}{^{\tau}}\right)  ^{k}}%
E_{k-p}\left(  w\right)  dw
\]
Further, we use the Beta function integral,
\[
\int_{0}^{\infty}\frac{w^{j}dw}{\left(  w+\frac{1}{^{\tau}}\right)  ^{k}}%
=\tau^{k}\int_{0}^{\infty}\frac{w^{j}dw}{\left(  \tau w+1\right)  ^{k}}%
=\tau^{k-j-1}\int_{0}^{\infty}\frac{u^{j}du}{\left(  u+1\right)  ^{k}}%
=\tau^{k-j-1}\frac{\Gamma\left(  j+1\right)  \Gamma\left(  k-j-1\right)
}{\Gamma\left(  k\right)  }%
\]
so that%
\[
L_{k}\left(  \tau\right)  =\sum_{j=0}^{p-1}\frac{\left(  -1\right)
^{j}\left(  \left(  k-j-2\right)  !\right)  ^{2}j!}{\left(  \left(
k-1\right)  !\right)  ^{2}}\tau^{k-j-1}+\frac{\left(  -1\right)  ^{p}\left(
k-p-1\right)  !}{\left(  k-1\right)  !}\int_{0}^{\infty}\frac{w^{p}e^{w}%
}{\left(  w+\frac{1}{^{\tau}}\right)  ^{k}}E_{k-p}\left(  w\right)  dw
\]
Finally, if $p=k-1$, then%
\[
L_{k}\left(  \tau\right)  =\sum_{j=0}^{k-2}\frac{\left(  -1\right)
^{j}\left(  \left(  k-j-2\right)  !\right)  ^{2}j!}{\left(  \left(
k-1\right)  !\right)  ^{2}}\tau^{k-j-1}+\frac{\left(  -1\right)  ^{k-1}%
}{\left(  k-1\right)  !}\int_{0}^{\infty}\frac{w^{k-1}e^{w}}{\left(
w+\frac{1}{^{\tau}}\right)  ^{k}}E_{1}\left(  w\right)  dw
\]
Introduced in (\ref{dominant_betaF_integral_expanded_Lk(tau)}) yields%
\begin{align*}
\int_{0}^{1}\left(  \int_{0}^{\infty}\frac{\left(  yu+1\right)  ^{N}%
e^{-\frac{1}{\tau}u}}{yu+1-y}du\right)  dy  &  =N!\sum_{k=0}^{N}\frac
{1}{\left(  N-k\right)  !}\sum_{j=0}^{k-1}\frac{\left(  -1\right)  ^{j}\left(
\left(  k-j-1\right)  !\right)  ^{2}j!}{\left(  k!\right)  ^{2}}\tau^{k-j}\\
&  +\sum_{k=0}^{N}\frac{N!}{\left(  N-k\right)  !}\frac{\left(  -1\right)
^{k}}{k!}\int_{0}^{\infty}\frac{w^{k}e^{w}}{\left(  w+\frac{1}{^{\tau}%
}\right)  ^{k+1}}E_{1}\left(  w\right)  dw
\end{align*}
The first sum can be rewritten as%
\begin{align*}
T  &  =\sum_{k=1}^{N}\frac{1}{\left(  N-k\right)  !}\sum_{j=0}^{k-1}%
\frac{\left(  -1\right)  ^{j}\left(  \left(  k-j-1\right)  !\right)  ^{2}%
j!}{\left(  k!\right)  ^{2}}\tau^{k-j}\\
&  =\sum_{k=1}^{N}\frac{1}{\left(  N-k\right)  !}\sum_{m=1}^{k}\frac{\left(
-1\right)  ^{k-m}\left(  \left(  m-1\right)  !\right)  ^{2}\left(  k-m\right)
!}{\left(  k!\right)  ^{2}}\tau^{m}\\
&  =\sum_{m=1}^{N}\left(  \sum_{k=m}^{N}\frac{1}{\left(  N-k\right)  !}%
\frac{\left(  -1\right)  ^{k-m}\left(  k-m\right)  !}{\left(  k!\right)  ^{2}%
}\right)  \left(  \left(  m-1\right)  !\right)  ^{2}\tau^{m}\\
&  =\frac{1}{N!}\sum_{m=1}^{N}\left(  \sum_{k=m}^{N}\binom{N}{k}\frac{\left(
-1\right)  ^{k-m}\left(  k-m\right)  !}{k!}\right)  \left(  \left(
m-1\right)  !\right)  ^{2}\tau^{m}%
\end{align*}
while the last sum equals%
\begin{align*}
\sum_{k=0}^{N}\frac{N!}{\left(  N-k\right)  !}\frac{\left(  -1\right)  ^{k}%
}{k!}\int_{0}^{\infty}\frac{w^{k}e^{w}}{\left(  w+\frac{1}{^{\tau}}\right)
^{k+1}}E_{1}\left(  w\right)  dw  &  =\int_{0}^{\infty}\frac{e^{w}E_{1}\left(
w\right)  }{w+\frac{1}{^{\tau}}}\sum_{k=0}^{N}\binom{N}{k}\left(  \frac
{-w}{w+\frac{1}{^{\tau}}}\right)  ^{k}dw\\
&  =\frac{1}{^{\tau^{N}}}\int_{0}^{\infty}\frac{e^{w}E_{1}\left(  w\right)
}{\left(  w+\frac{1}{^{\tau}}\right)  ^{N+1}}dw
\end{align*}
After combining the above, we arrive at%
\begin{equation}
\int_{0}^{1}\left(  \int_{0}^{\infty}\frac{\left(  yu+1\right)  ^{N}%
e^{-\frac{1}{\tau}u}}{yu+1-y}du\right)  dy=\sum_{m=1}^{N}C_{m}\tau^{m}%
+\tau\int_{0}^{\infty}\frac{e^{w}E_{1}\left(  w\right)  }{\left(  \tau
w+1\right)  ^{N+1}}dw \label{beta_F(tau)_integral_Taylorseries_Cj}%
\end{equation}
with%
\begin{equation}
C_{m}=\left(  \left(  m-1\right)  !\right)  ^{2}\sum_{k=m}^{N}\binom{N}%
{k}\frac{\left(  -1\right)  ^{k-m}\left(  k-m\right)  !}{k!}
\label{def_coeff_Cm}%
\end{equation}

\begin{theorem}
\label{theorem_Taylor_series_Bj_Cj} The Taylor series of $\beta F\left(
\tau\right)  $ is defined by%
\begin{equation}
\beta F\left(  \tau\right)  =\sum_{j=1}^{N}B_{j}\tau^{j}
\label{beta_F(tau)_Taylor_series_coefficients_Bj}%
\end{equation}
with%
\begin{equation}
B_{j}=\sum_{k=j}^{N}\frac{\left(  N-k+j-1\right)  !}{\left(  N-k\right)  !k}
\label{def_coeff_Bj}%
\end{equation}
Moreover, another form for $B_{j}$ is%
\begin{equation}
B_{j}=\left(  \left(  j-1\right)  !\right)  ^{2}\sum_{k=j}^{N}\binom{N}%
{k}\frac{\left(  -1\right)  ^{k-j}\left(  k-j\right)  !}{k!} \label{Bj=Cj}%
\end{equation}

\end{theorem}

\textbf{Proof:} Reversing the $j$- and $r$-sum in (\ref{def_beta_F}) yields%
\[
\beta F\left(  \tau\right)  =\sum_{j=1}^{N}\left(  \sum_{r=j}^{N}\frac{\left(
N-r+j-1\right)  !}{\left(  N-r\right)  !r}\right)  \tau^{j}%
\]
which proves (\ref{beta_F(tau)_Taylor_series_coefficients_Bj}) and
(\ref{def_coeff_Bj}). We will now prove that $B_{j}=C_{j}$, where $C_{j}$ is
defined in (\ref{def_coeff_Cm}). First, we prove that%
\[
\int_{0}^{\infty}\frac{e^{w}E_{n}\left(  w\right)  }{\left(  \tau w+1\right)
}dw=\int_{0}^{\infty}\frac{e^{w}E_{1}\left(  w\right)  }{\left(  \tau
w+1\right)  ^{n}}dw
\]
Introducing the definition of the exponential integral of order $n$ yields%
\begin{align*}
\int_{0}^{\infty}\frac{e^{w}E_{n}\left(  w\right)  }{\left(  \tau w+1\right)
}dw  &  =\int_{0}^{\infty}\frac{e^{w}dw}{\left(  \tau w+1\right)  }\int
_{1}^{\infty}\frac{e^{-wt}}{t^{n}}dt\\
&  =\int_{1}^{\infty}\frac{dt}{t^{n}}\int_{0}^{\infty}\frac{e^{-w\left(
t-1\right)  }dw}{\left(  \tau w+1\right)  }\\
&  =\int_{0}^{\infty}\frac{du}{\left(  u+1\right)  ^{n}}\int_{0}^{\infty}%
\frac{e^{-wu}dw}{\left(  \tau w+1\right)  }%
\end{align*}
Let $y=\tau w+1$, then%
\begin{align*}
\int_{0}^{\infty}\frac{e^{w}E_{n}\left(  w\right)  }{\left(  \tau w+1\right)
}dw  &  =\frac{1}{\tau}\int_{0}^{\infty}\frac{du}{\left(  u+1\right)  ^{n}%
}\int_{1}^{\infty}\frac{e^{-\left(  y-1\right)  \frac{u}{\tau}}dy}{y}\\
&  =\frac{1}{\tau}\int_{0}^{\infty}\frac{e^{\frac{u}{\tau}}du}{\left(
u+1\right)  ^{n}}\int_{1}^{\infty}\frac{e^{-y\frac{u}{\tau}}dy}{y}\\
&  =\frac{1}{\tau}\int_{0}^{\infty}\frac{e^{\frac{u}{\tau}}E_{1}\left(
\frac{u}{\tau}\right)  }{\left(  u+1\right)  ^{n}}du=\int_{0}^{\infty}%
\frac{e^{w}E_{1}\left(  w\right)  }{\left(  \tau w+1\right)  ^{n}}du
\end{align*}

From (\ref{tweede_stuk_beta_F(tau)_intergral}), we find that%
\[
\beta F\left(  \tau\right)  =\int_{0}^{1}\left(  \int_{0}^{\infty}%
\frac{\left(  yu+1\right)  ^{N}e^{-\frac{1}{\tau}u}}{yu+1-y}du\right)
dy-\tau\int_{0}^{\infty}\frac{e^{w}E_{1}\left(  w\right)  }{\left(  \tau
w+1\right)  ^{N+1}}dw
\]
The uniqueness of the Taylor series implies that
(\ref{beta_F(tau)_Taylor_series_coefficients_Bj}) equals the Taylor series in
(\ref{beta_F(tau)_integral_Taylorseries_Cj}).\hfill$\square\medskip$

\subsubsection{The coefficients $B_{j}$}

We present other properties of the Taylor coefficients $B_{j}$, defined in
(\ref{def_coeff_Bj}). Starting from%
\[
\sum_{r=1}^{N}\frac{\left(  N+j-1-r\right)  !}{\left(  N-r\right)  !r}%
=\sum_{r=1}^{j-1}\frac{\left(  N+j-1-r\right)  !}{\left(  N-r\right)  !r}%
+\sum_{r=j}^{N}\frac{\left(  N+j-1-r\right)  !}{\left(  N-r\right)  !r}%
\]
and using (see \cite[p. 410]{PVM_PerformanceAnalysisCUP}),
\[
\sum_{r=1}^{N}\frac{\left(  N+j-1-r\right)  !}{\left(  N-r\right)  !r}%
=\frac{\left(  N+j-1\right)  !}{N!}\left[  \psi(N+j)-\psi(j)\right]
=\frac{\left(  N+j-1\right)  !}{N!}\sum_{k=j}^{N+j-1}\frac{1}{k}%
\]
shows that%
\begin{equation}
B_{j}=\frac{\left(  N+j-1\right)  !}{N!}\sum_{k=j}^{N+j-1}\frac{1}{k}%
-\sum_{r=1}^{j-1}\frac{\left(  N+j-1-r\right)  !}{\left(  N-r\right)  !r}
\label{coeff_Bj_1}%
\end{equation}

A recursion for the Taylor coefficient%
\[
B_{j}\left(  N\right)  =\sum_{r=j}^{N}\frac{\left(  N+j-1-r\right)  !}{\left(
N-r\right)  !r}%
\]
is%
\begin{equation}
B_{j+1}\left(  N\right)  =B_{j+1}\left(  N-1\right)  +jB_{j}\left(  N\right)
-\frac{\left(  N-1\right)  !}{\left(  N-j\right)  !} \label{recursion_Bj(N)}%
\end{equation}
Indeed,%
\begin{align*}
B_{j+1}\left(  N\right)   &  =\sum_{r=j+1}^{N}\frac{\left(  N+j-r\right)
!}{\left(  N-r\right)  !r}=\sum_{r=j+1}^{N}\frac{\left(  N+j-r\right)  \left(
N+j-1-r\right)  !}{\left(  N-r\right)  !r}\\
&  =\sum_{r=j+1}^{N}\frac{\left(  N+j-1-r\right)  !}{\left(  N-r-1\right)
!r}+j\sum_{r=j+1}^{N}\frac{\left(  N+j-1-r\right)  !}{\left(  N-r\right)
!r}\\
&  =B_{j+1}\left(  N-1\right)  +j\left(  B_{j}\left(  N\right)  -\frac{\left(
N-1\right)  !}{\left(  N-j\right)  !j}\right)
\end{align*}

\subsection{Order estimate of $F\left(  \tau\right)  $ for large $N$}

\label{sec_order_F(tau_largeN}The exponential integral is bounded
\cite[5.1.9]{Abramowitz} by%
\begin{equation}
\frac{1}{x+n}<e^{x}E_{n}\left(  x\right)  \leq\frac{1}{x+n-1}
\label{bounds_exponential_integral}%
\end{equation}
The bound (\ref{bounds_exponential_integral}) will play a crucial role in the
proof of Theorem \ref{theorem_bound_f1_eps=0} below. In the determination of
the order of $F\left(  \tau\right)  $ for large $N$ (and fixed $x=N\tau$), the
bound (\ref{bounds_exponential_integral}) can be used when the last term in
the summation is treated separately
\begin{align*}
\int_{0}^{1}\left(  \int_{0}^{\infty}\frac{\left(  yu+1\right)  ^{N}%
e^{-\frac{1}{\tau}u}}{yu+1-y}du\right)  dy  &  =N!\sum_{j=0}^{N}\frac{1}%
{j!}\int_{0}^{\infty}\frac{e^{w}E_{N+1-j}\left(  w\right)  }{\left(
w+\frac{1}{^{\tau}}\right)  ^{N+1-j}}dw\\
&  =N!\sum_{j=0}^{N-1}\frac{1}{j!}\int_{0}^{\infty}\frac{e^{w}E_{N+1-j}\left(
w\right)  }{\left(  w+\frac{1}{^{\tau}}\right)  ^{N+1-j}}dw+\int_{0}^{\infty
}\frac{e^{w}E_{1}\left(  w\right)  }{w+\frac{1}{^{\tau}}}dw
\end{align*}
The definition (\ref{def_Lk(tau)}) of $L_{k}\left(  \tau\right)  $ and the
inequality $E_{k}\left(  w\right)  <E_{k-1}\left(  w\right)  $ shows that, for
$k>1$,%
\[
L_{k}\left(  \tau\right)  =\int_{0}^{\infty}\frac{e^{w}E_{k}\left(  w\right)
}{\left(  w+\frac{1}{^{\tau}}\right)  ^{k}}dw<\int_{0}^{\infty}\frac
{e^{w}E_{k-1}\left(  w\right)  }{\left(  w+\frac{1}{^{\tau}}\right)  ^{k-1}%
}\frac{1}{\left(  w+\frac{1}{^{\tau}}\right)  }dw<\tau L_{k-1}\left(
\tau\right)
\]
and%
\[
\frac{dL_{k}\left(  \tau\right)  }{d\tau}=\frac{k}{\tau^{2}}\int_{0}^{\infty
}\frac{e^{w}E_{k}\left(  w\right)  }{\left(  w+\frac{1}{^{\tau}}\right)
^{k+1}}dw>\frac{k}{\tau^{2}}L_{k+1}\left(  \tau\right)
\]
Using $\frac{dE_{n}\left(  x\right)  }{dx}=-E_{n-1}\left(  x\right)  $, after
partial integration of $L_{k}\left(  \tau\right)  $, we obtain, for $k>1$,%
\[
L_{k}\left(  \tau\right)  =\frac{\tau^{k-1}}{\left(  k-1\right)  ^{2}}%
+\frac{1}{\left(  k-1\right)  }\int_{0}^{\infty}\frac{e^{w}\left(
E_{k}\left(  w\right)  -E_{k-1}\left(  w\right)  \right)  }{\left(  w+\frac
{1}{^{\tau}}\right)  ^{k-1}}dw
\]
Since $E_{k}\left(  w\right)  <E_{k-1}\left(  w\right)  $, we find, for $k>1$,
that
\begin{equation}
L_{k}\left(  \tau\right)  <\frac{\tau^{k-1}}{\left(  k-1\right)  ^{2}}
\label{upper_bound_L_k}%
\end{equation}
The function $L_{1}\left(  \tau\right)  $ requires a different treatment,%
\begin{align*}
L_{1}\left(  \tau\right)   &  =\int_{0}^{\infty}\frac{e^{w}E_{1}\left(
w\right)  }{w+\frac{1}{^{\tau}}}dw=\int_{0}^{\infty}\frac{e^{w}dw}{\left(
w+\frac{1}{^{\tau}}\right)  }\int_{w}^{\infty}\frac{e^{-t}}{t}dt\\
&  =\int_{0}^{\infty}dw\int_{w}^{\infty}dt\frac{e^{-\left(  t-w\right)  }%
}{\left(  w+\frac{1}{^{\tau}}\right)  t}=\int_{0}^{\infty}dw\int_{0}^{\infty
}dy\frac{e^{-y}}{\left(  w+\frac{1}{^{\tau}}\right)  \left(  y+w\right)  }\\
&  =\int_{0}^{\infty}e^{-y}dy\int_{0}^{\infty}\frac{dw}{\left(  w+\frac
{1}{^{\tau}}\right)  \left(  y+w\right)  }%
\end{align*}
The integral, computed after partial fraction expansion,%
\begin{align*}
\int_{0}^{\infty}\frac{dw}{\left(  w+\frac{1}{^{\tau}}\right)  \left(
y+w\right)  }  &  =\frac{1}{\left(  \frac{1}{^{\tau}}-y\right)  }\left\{
\int_{0}^{\infty}\frac{dw}{\left(  y+w\right)  }-\int_{0}^{\infty}\frac
{dw}{\left(  w+\frac{1}{^{\tau}}\right)  }\right\} \\
&  =\frac{\ln\left(  y\tau\right)  }{y-\frac{1}{^{\tau}}}=\frac{\ln\left(
y\right)  -\ln\left(  \frac{1}{^{\tau}}\right)  }{y-\frac{1}{^{\tau}}}%
\end{align*}
shows that $\frac{\ln\left(  y\tau\right)  }{y-\frac{1}{^{\tau}}}$ is
decreasing in $y$ and $\frac{1}{^{\tau}}$. Moreover, we find that%
\[
L_{1}\left(  \tau\right)  =\int_{0}^{\infty}\frac{\ln\left(  y\tau\right)
e^{-y}}{y-\frac{1}{^{\tau}}}dy=\int_{0}^{\infty}\frac{\ln\left(  u\right)
e^{-\frac{u}{\tau}}}{u-1}du
\]
Thus,%
\[
L_{1}\left(  \tau\right)  =\int_{0}^{\infty}\frac{e^{w}E_{1}\left(  w\right)
}{\left(  w+\frac{1}{^{\tau}}\right)  }dw=\int_{0}^{\infty}\frac{\ln\left(
u\right)  e^{-\frac{u}{\tau}}}{u-1}du
\]
Now, let $t=\frac{1}{\tau}$, then%
\begin{align*}
\frac{dL_{1}\left(  t\right)  }{dt}  &  =-\int_{0}^{\infty}\frac{u\ln\left(
u\right)  e^{-tu}}{u-1}du=-\int_{0}^{\infty}\frac{\left(  u-1+1\right)
\ln\left(  u\right)  e^{-tu}}{u-1}du\\
&  =-\int_{0}^{\infty}\ln\left(  u\right)  e^{-tu}du-\int_{0}^{\infty}%
\frac{\ln\left(  u\right)  e^{-tu}}{u-1}du
\end{align*}
With%
\begin{align*}
\int_{0}^{\infty}\ln\left(  u\right)  e^{-tu}du  &  =\frac{1}{t}\int
_{0}^{\infty}\ln\left(  \frac{y}{t}\right)  e^{-y}dy=\frac{1}{t}\int
_{0}^{\infty}\ln\left(  y\right)  e^{-y}dy-\frac{\ln t}{t}\int_{0}^{\infty
}e^{-y}dy\\
&  =-\frac{\gamma}{t}-\frac{\ln t}{t}%
\end{align*}
because $\int_{0}^{\infty}e^{-y}\ln ydy=\Gamma^{\prime}\left(  1\right)
=-\gamma=-0.5227$, we obtain the first-order differential equation%
\[
\frac{dL_{1}\left(  t\right)  }{dt}=\frac{\gamma+\ln t}{t}-L_{1}\left(
t\right)
\]
Since $\frac{u\ln\left(  u\right)  }{u-1}\leq\sqrt{u}$ (with equality for
$u=1$ and tight for $0\leq u\leq1$),%
\[
0\leq-\frac{dL_{1}\left(  t\right)  }{dt}<\int_{0}^{\infty}\sqrt{u}%
e^{-tu}du=\frac{1}{t\sqrt{t}}\Gamma\left(  \frac{3}{2}\right)
\]
the differential equation provides us with%
\[
\tau\left(  \gamma-\ln\tau\right)  <L_{1}\left(  \tau\right)  <\tau\left(
\gamma-\ln\tau\right)  +\tau^{\frac{3}{2}}\frac{\sqrt{\pi}}{2}%
\]
and the upper bound is tight for $\tau<1$, but loose for $\tau>1$. Since we
are interested in small $\tau=\frac{x}{N}$, the upper bound suffices and
illustrates that $L_{1}\left(  \frac{x}{N}\right)  =O\left(  \frac{\ln N}%
{N}\right)  $.

The main result here is the following theorem

\begin{theorem}
\label{theorem_bound_f1_eps=0} For $\tau=\frac{x}{N}$, $F\left(  \tau\right)
=\lim_{\varepsilon\rightarrow0}f_{1}$ behaves for large $N$ and fixed $x>1$ as%
\[
F\left(  \frac{x}{N}\right)  \sim\frac{1}{\delta}\frac{x\sqrt{2\pi}}{\left(
x-1\right)  ^{2}}\frac{\exp\left(  N\left\{  \log x+\frac{1}{x}-1\right\}
\right)  }{\sqrt{N}}%
\]

\end{theorem}

\textbf{Proof:} Using the bounds (\ref{bounds_exponential_integral}) for the
exponential integral, the negative term in
(\ref{betaF_in_exponential_integrals})
\[
\int_{0}^{\infty}\frac{dw}{\left(  w+\frac{1}{\tau}\right)  \left(
w+N-1\right)  }<\int_{0}^{\infty}\frac{e^{w}E_{N+1}\left(  w\right)
dw}{w+\frac{1}{\tau}}\leq\int_{0}^{\infty}\frac{dw}{\left(  w+\frac{1}{\tau
}\right)  \left(  w+N\right)  }%
\]
With%
\begin{align*}
\int_{0}^{\infty}\frac{dw}{\left(  w+\frac{1}{\tau}\right)  \left(
w+m\right)  }  &  =-\left(  m-\frac{1}{\tau}\right)  \left.  \ln\frac
{w+m}{w+\frac{1}{\tau}}\right\vert _{0}^{\infty}\\
&  =\left(  m-\frac{1}{\tau}\right)  \ln m\tau
\end{align*}
we arrive at%
\[
\left(  N-1-\frac{1}{\tau}\right)  \ln\left(  N-1\right)  \tau<\int_{0}%
^{1}y^{N}\left(  \int_{0}^{\infty}\frac{e^{-\frac{1}{\tau}u}}{yu+1-y}%
du\right)  dy\leq\left(  N-\frac{1}{\tau}\right)  \ln N\tau
\]
which demonstrates, since $x=N\tau>1$ and fixed, that the negative term in
(\ref{betaF_in_exponential_integrals}) grows as $N\left(  1-\frac{1}%
{x}\right)  \log x$ for large $N$.

Using the bounds (\ref{bounds_exponential_integral}) in the expression
(\ref{def_Lk(tau)}) of $L_{k}\left(  \tau\right)  $ yields, for $0\leq j<N$,%
\[
\int_{0}^{\infty}\frac{dw}{\left(  w+\frac{1}{^{\tau}}\right)  ^{N+1-j}\left(
w+N+1-j\right)  }<\int_{0}^{\infty}\frac{e^{w}E_{N+1-j}\left(  w\right)
dw}{\left(  w+\frac{1}{^{\tau}}\right)  ^{N+1-j}}\leq\int_{0}^{\infty}%
\frac{dw}{\left(  w+\frac{1}{^{\tau}}\right)  ^{N+1-j}\left(  w+N-j\right)  }%
\]
The upper bound is further%
\begin{align*}
\int_{0}^{\infty}\frac{e^{w}E_{N+1-j}\left(  w\right)  dw}{\left(  w+\frac
{1}{^{\tau}}\right)  ^{N+1-j}}  &  \leq\int_{0}^{\infty}\frac{dw}{\left(
w+\frac{1}{^{\tau}}\right)  ^{N+1-j}\left(  w+N-j\right)  }<\frac{1}{N-j}%
\int_{0}^{\infty}\frac{dw}{\left(  w+\frac{1}{^{\tau}}\right)  ^{N+1-j}}\\
&  =\frac{\tau^{N-j}}{\left(  N-j\right)  ^{2}}%
\end{align*}
which also follows from (\ref{upper_bound_L_k}). Hence,%
\begin{equation}
\int_{0}^{1}\left(  \int_{0}^{\infty}\frac{\left(  yu+1\right)  ^{N}%
e^{-\frac{1}{\tau}u}}{yu+1-y}du\right)  dy<L_{1}\left(  \tau\right)
+N!\tau^{N}\sum_{j=0}^{N-1}\frac{\tau^{-j}}{j!}\frac{1}{\left(  N-j\right)
^{2}} \label{upper_bound_first_integral_betaF}%
\end{equation}
When making the rather crude approximation%
\[
\sum_{j=0}^{N-1}\frac{\tau^{-j}}{j!}\frac{1}{\left(  N-j\right)  ^{2}}%
<\sum_{j=0}^{N-1}\frac{\tau^{-j}}{j!}<e^{\frac{1}{\tau}}%
\]
we obtain the upper bound%
\begin{equation}
\beta F\left(  \tau\right)  <L_{1}\left(  \tau\right)  +N!\tau^{N}e^{\frac
{1}{\tau}} \label{bound_f1_eps=0}%
\end{equation}

A much better approximation, that is asymptotically correct for large $N$ and
constant $x>1$, is%
\[
\sum_{j=0}^{N-1}\frac{\left(  \frac{N}{x}\right)  ^{j}}{j!}\frac{1}{\left(
N-j\right)  ^{2}}\approx\frac{e^{\frac{N}{x}}}{\left(  1-\frac{1}{x}\right)
^{2}N^{2}}%
\]
The arguments require order estimates for sums $\sum_{m=a}^{b}u_{m}$, where
$u_{m}=e^{-y}\frac{y^{m}}{m!}$, derived in a theorem by Hardy \cite[p.
200]{Hardy_div}. Since $\tau=\frac{x}{N}$, with $x=1+\delta>1$, we obtain%
\[
\sum_{j=0}^{N-1}\frac{\tau^{-j}}{j!}=\sum_{j=0}^{N-1}\frac{\left(  \frac{N}%
{x}\right)  ^{j}}{j!}=\sum_{j=0}^{\frac{N}{x}\left(  1+\delta\right)  -1}%
\frac{\left(  \frac{N}{x}\right)  ^{j}}{j!}%
\]
The largest term in the $j$-sum occurs at $j=\frac{N}{x}$. Hardy shows, that
for large $N$ and $0<\delta<1$, the sum is very close to $e^{\frac{N}{x}}$
with error at most $O\left(  e^{-\frac{1}{3}\delta^{2}}\right)  $. Moreover,
for $\left\vert h\right\vert \leq y^{\zeta}$ (and $\frac{1}{2}<\zeta<\frac
{2}{3}$), $u_{m}=e^{-y}\frac{y^{m}}{m!}$ with $m=\left[  y\right]  +h$ equals%
\[
u_{m}=\sqrt{\frac{1}{2\pi\left[  y\right]  }}e^{-\frac{h^{2}}{2\left[
y\right]  }}\left\{  1+O\left(  \frac{\left\vert h\right\vert +1}{y}\right)
+O\left(  \frac{h^{3}}{y^{2}}\right)  \right\}
\]
Since the terms $u_{\left[  y\right]  +h}$ are increasingly peaked around the
maximum $y=\frac{N}{x}$, we can approximate%
\[
e^{-\frac{1}{\tau}}\sum_{j=0}^{N-1}\frac{\tau^{-j}}{j!}\frac{1}{\left(
N-j\right)  ^{2}}\approx\frac{1}{\left(  N-\left[  \frac{N}{x}\right]
\right)  ^{2}}\sum_{h=-\left[  \frac{N}{x}\right]  ^{\zeta}}^{\left[  \frac
{N}{x}\right]  ^{\zeta}}\frac{e^{-\frac{xh^{2}}{2N}}}{\sqrt{2\pi\left[
\frac{N}{x}\right]  }}\approx\frac{1}{N^{2}\left(  1-\frac{1}{x}\right)  ^{2}}%
\]
Hence, for large $N$ and fixed $x=\tau N$ in the dominant term
(\ref{upper_bound_first_integral_betaF}) for $\beta F\left(  \tau\right)  $,
we arrive at%
\begin{align*}
\beta F\left(  \frac{N}{x}\right)   &  \approx L_{1}\left(  \frac{x}%
{N}\right)  +N!\left(  \frac{x}{N}\right)  ^{N}\sum_{j=0}^{N-1}\frac{\left(
\frac{N}{x}\right)  ^{j}}{j!}\frac{1}{\left(  N-j\right)  ^{2}}\\
&  \sim O\left(  \frac{\ln N}{N}\right)  +N!\left(  \frac{x}{N}\right)
^{N}\frac{e^{\frac{N}{x}}}{\left(  1-\frac{1}{x}\right)  ^{2}N^{2}}%
\end{align*}
Using Stirling's approximation $N!\approx\sqrt{2\pi N}N^{N}e^{-N}$ and
$\beta=\tau\delta=\frac{x}{N}\delta$, we finally obtain
\[
F\left(  \frac{N}{x}\right)  \sim\frac{1}{\delta}\frac{N}{x}\frac{\sqrt{2\pi
N}N^{N}e^{-N}}{N^{N}\left(  1-\frac{1}{x}\right)  ^{2}N^{2}}\left(
xe^{\frac{1}{x}}\right)  ^{N}%
\]
which proves the theorem after some simplifications.\hfill$\square\medskip$

\end{document}